\documentclass[a4paper,12pt,reqno]{amsart}
\usepackage{mathrsfs}
\usepackage{amsmath,amssymb,latexsym,amsfonts,amscd}

\title{Explicit birational geometry of threefolds of general type, I}
\author{Jungkai A. Chen and Meng Chen}
\address{\rm Department of Mathematics, National Taiwan University, Taipei,
106, Taiwan} \email{jkchen@math.ntu.edu.tw}

\address{\rm Institute of Mathematics, Fudan University,
Shanghai, 200433, People's Republic of China}
\email{mchen@fudan.edu.cn}

\thanks{The first author was partially supported by TIMS, NCTS/TPE
and National Science Council of Taiwan. The second author was
supported by National Outstanding Young Scientist Foundation
(\#10625103) and NNSFC Key project (\#10731030)}
\numberwithin{equation}{section}

\newcommand{\bC}{{\mathbb C}}
\newcommand{\bQ}{{\mathbb Q}}
\newcommand{\bP}{{\mathbb P}}

\newcommand\Vol{\text{\rm Vol}}

\newcommand\OO{{\mathcal{O}}}

\newtheorem{thm}{Theorem}[section]
\newtheorem{lem}[thm]{Lemma}
\newtheorem{cor}[thm]{Corollary}

\theoremstyle{definition}
\newtheorem{defn}[thm]{Definition}
\newtheorem{setup}[thm]{}

\newtheorem{exmp}[thm]{Example}

\newtheorem{rem}[thm]{Remark}
\theoremstyle{remark}

\newtheorem{Prbm}{\bf Problem}

\begin{document}
\begin{abstract}
Let $V$ be a complex nonsingular projective 3-fold of general type.
We prove $P_{12}(V):=\dim H^0(V, 12K_V)>0$ and $P_{m_0}(V)>1$ for
some positive integer $m_0\leq 24$. A direct consequence is the
birationality of the pluricanonical map $\varphi_m$ for all $m\geq
126$. Besides, the canonical volume $\text{Vol}(V)$ has a universal
lower bound $\nu(3)\geq \frac{1}{63\cdot 126^2}$.
\end{abstract}

\maketitle
\pagestyle{myheadings} \markboth{\hfill J. A. Chen and M. Chen
\hfill}{\hfill Explicit birational geometry of threefolds\hfill}




{\sf R\'esum\'e en fran\c{c}ais:} Soit $V$ une vari\'et\'e non singuli\`ere complexe
de type g\'en\'eral  de dimension 3. Nous montrons   $P_{12}(V):=\dim H^0(V, 12K_V)>0$ et $P_{m_0}(V)>1$ pour un certain entier $m_0\leq 24$. Une cons\'equence directe est la birationalit\'e de l'application pluricanonique
$\varphi_m$ pour tout $m\geq
126$. De plus, le volume canonique $\text{Vol}(V)$ a un minorant universel $\nu(3)\geq \frac{1}{63\cdot 126^2}$. \\[.2cm]

\section{\bf Introduction}
Let $Y$ be a nonsingular projective variety of dimension $n$. It is
said to be of general type if the pluricanonical map
$\varphi_m:=\Phi_{|mK_Y|}$ corresponding to the linear system
$|mK_Y|$ is birational into a projective space for $m \gg 0$. Thus
it is natural and important to ask:

\begin{Prbm}\label{Pb} Can one find a constant $c(n)$, so that $\varphi_{m}$ is
birational onto its image for all $m \ge c(n)$ and for all $Y$ with
$\dim Y =n$?
\end{Prbm}

When $\dim Y=1$, it was classically known that $|mK_Y|$ gives an
embedding of $Y$ into a projective space if $m \ge 3$. When $\dim
Y=2$, Kodaira-Bombieri's theorem \cite{Bom} says that $|mK_Y|$ gives
a birational map onto the image for $m\geq 5$. This theorem has
essentially established the canonical classification theory for
surfaces of general type.

A natural approach to study this problem in higher dimensions is an
induction on the dimension by utilizing vanishing theorems. This
amounts to estimating the plurigenus, for which purpose the greatest
difficulty seems to be to bound from below the canonical volume
$$\Vol (Y):=\limsup_{\{m\in {\mathbb Z}^{+}\}}
\{\frac{n!}{m^n}\dim_{\bC}H^0(Y, \OO_Y(mK_Y))\}.$$ The volume is an
integer when $\dim Y \le 2$. However it's only a rational number in
general, which may account for the complexity of high dimensional
birational geometry. In fact, it is almost an equivalent question
to study the lower bound of the canonical volume.

\begin{Prbm} Can one find a constant $\nu(n)$ such that $\Vol(Y) \geq
\nu(n)$ for all varieties $Y$ of general type with $\dim Y =n$?
\end{Prbm}

A recent result of Hacon and M$^{\text{c}}$Kernan \cite{H-M},
Takayama \cite{Ta} and Tsuji \cite{Tsuji} shows the existence of
both $c(n)$ and $\nu(n)$. An explicit constant $c(n)$ or $ \nu (n)$
is, however, mysterious at least up to now. Notice that similar
questions were asked by Koll\'ar and Mori \cite[7.74]{K-M}.

Here we mainly deal with $c(3)$ and $\nu(3)$. For known results
under extra assumptions, one may refer to \cite{Crelle, JC-H, JLMS,
IJM, MathAnn, Chen-Zuo, Ein-Lazars, C-R, Kol, Luo} and others.  In
this series of papers, we would like to present two realistic
constants $c(3)$ and $\nu(3)$. In fact, our method can help us to
prove some sharp results. Being worried about that a very long
paper will cause the tiredness to readers, we decide to only
explain our key technique and rough statements in the first part
whereas more refined and some sharp statements will be presented in
the subsequent papers. Our main result in this paper is the
following:

\begin{thm}\label{main} Let $V$ be a nonsingular projective
3-fold of general type. Then
\begin{itemize}
\item[(1)] $P_{12}>0$;

\item[(2)] $P_{m_0}\geq 2$ for some positive integer $m_0\leq 24$.
\end{itemize}
\end{thm}

With Koll\'ar's result \cite[Corollary 4.8]{Kol} and its improved
form \cite[Theorem 0.1]{JPAA}, we immediately get the following:

\begin{cor}\label{birat} Let $V$ be a nonsingular projective 3-fold of general
type. Then
\begin{itemize}
\item[(1)] $\varphi_m$ is birational onto its image for all $m\geq
126$.

\item[(2)] $\text{Vol}(V)\geq \frac{1}{63\cdot 126^2}$.
\end{itemize}
\end{cor}

\begin{exmp} (see \cite[p151, No.23]{C-R}) The ``worst'' known example is a general weighted
hypersurface $X=X_{46} \subset \mathbb{P}(4,5,6,7,23)$. The 3-fold
$X$ has invariants: $p_g(X)=P_2(X)=P_3(X)=0$,
$P_4(X)=\cdots=P_9(X)=1$, $P_{10}(X)=2$ and
$\Vol(X)=\frac{1}{420}$. Moreover, it is known that $\varphi_{m}$
is birational for all $m\geq 27$, but $\varphi_{26}$ is not
birational.
\end{exmp}

Now we explain the main idea of our paper. It is very natural to
investigate the plurigenus $P_m$, which can be calculated using
Reid's Riemann-Roch formula in \cite{C-3f, YPG}. However the most
difficult point is to control the contribution from singularities
due to the combinatorial complexity of baskets of singularities on
the 3-fold.

Indeed, given a minimal 3-fold $X$ with at worst canonical
singularities, a known fact is that the canonical volume and all
plurigenera are determined by the basket (of singularities) $B$,
$\chi=\chi(\mathcal{O}_X)$ and $P_2=P_2(X)$. We call the triple
$(B,\chi, P_2)$ a {\it formal basket}. {}First we will define a
partial ordering (called ``packing'') between formal baskets. (In
this paper, we are only concerned about the numerical behavior of
``packing'', rather than its geometric meaning. More details on its
geometric aspect will be explored in our subsequent works.) Then we
introduce the ``canonical sequence of prime unpackings of a basket"
$$ B^{(0)} \succcurlyeq B^{(5)} \succcurlyeq ... \succcurlyeq B^{(n)} \succcurlyeq ... \succcurlyeq B$$
and, furthermore, each step in the sequence can be calculated in
terms of the datum of the given formal basket. The intrinsic
properties of the canonical sequence tell us many new inequalities
among the Euler characteristic and the plurigenus, of which the
most interesting one is:
$$2P_5+3P_6+P_8+P_{10}+P_{12} \ge \chi +10
P_2+4P_3+P_7+P_{11}+P_{13}.$$ If $P_{m_0}\geq 2$ for some $m_0\leq
12$, then one gets many interesting results by  \cite[Corollary
4.8]{Kol} and \cite[Theorem 0.1]{JPAA}. Otherwise one has $P_m\leq
1$ for all $m\leq 12$ and the above inequality gives $\chi\leq 8$.
This essentially tells us that the number of formal baskets is
finite! Thus, theoretically, we are able to obtain various effective
results.

Here is the overview to the structure of this paper. In Section 2,
we introduce the notion of packing and define some invariants of
baskets. Then we define the canonical sequence of ``prime
unpackings'' of a basket and give some examples. In Section 3, we
define the notion of formal baskets. Then we study various relations
among formal baskets, Euler characteristics and $K^3$. We calculate
the first few elements in the canonical sequence of the given
basket. This immediately gives many inequalities among Euler
characteristics. We would like to remark that the method so far
works for $\mathbb{Q}$-factorial threefolds (not only of general
type) with canonical singularities. With all these preparations, we
prove the main theorem on threefolds of general type in Section 4.

Another remark is that the method in Sections 2 and 3 is also valid
for $\bQ$-Fano threefolds. More precisely, there are similar
relations among formal baskets, anti-plurigenera and the
anti-canonical volume with proper sign alterations because of Serre
dualities. We will explore some more applications of our method in a
future work.

In our next paper of this series, we will work out some
classification of formal baskets with given small Euler
characteristics. Together with some more detailed study of the
geometry of pluricanonical maps, we will prove the following
theorem:
\medskip

\noindent{\bf Theorem A.} Let $V$ be a nonsingular projective 3-fold
of general type. Then the following holds.
\begin{itemize}
\item[(i)] $\varphi_m$ is birational onto its image for all $m\geq
73$.

\item[(ii)] $\Vol(V)\geq \frac{1}{2660}$.

\item[(iii)] Suppose that $\chi(\OO_V)\leq 1$. Then $\Vol(V)\geq \frac{1}{420}$, which
is optimal. Moreover $\varphi_m$ is birational for all $m\geq 40$.
\end{itemize}
\medskip

Throughout, we work over the complex number field $\bC$. We prefer
to use $\sim$ to denote the linear equivalence and $\equiv$ means
numerical equivalence. We mainly refer to \cite{KMM, K-M, Reid83}
for tool books on 3-dimensional birational geometry.
\medskip

\noindent {\bf Acknowledgments.} We would like to thank Gavin
Brown, Hou-Yi Chen, Jiun-Cheng Chen, H\'el\`ene Esnault, Christopher Hacon, J\'anos
Koll\'ar, Hui-Wen Lin, Miles Reid, Pei-Yu Tsai, Chin-Lung Wang and
De-Qi Zhang for their generous helps and comments on this subject.

\section{\bf Baskets of singularities}

In this section, we introduce the notion of packing  between baskets
of singularities. This notion defines a partial ordering on the set
of baskets.  For a given basket, we define its canonical sequence of
prime unpackings. The canonical sequence trick is a fundamental and
effective tool in our arguments.

\begin{setup}{\bf Terminal quotient singularity and basket.} By {\it a 3-dimensional terminal quotient singularity $Q$ of
type $\frac{1}{r}(1,-1,b)$}, we mean a singularity which is
analytically isomorphic to the quotient of $(\bC^3, \text{o})$ by a
cyclic group action  $\varepsilon$:
$$\varepsilon(x,y,z)=(\varepsilon x,\varepsilon^{-1} y, \varepsilon^b z)$$
where $r$ is a positive integer, $\varepsilon$ is a fixed $r$-th
primitive root of $1$, the integer $b$ is coprime to $r$ and
$0<b<r$.
\end{setup}

\begin{setup}{\bf Convention.} By replacing $\varepsilon$ with
another primitive root of 1 and changing the ordering of
coordinates, we may and will assume that $b \leq \frac{r}{2}$.
\end{setup}
{\it A basket} $B$ of singularities is a collection (allowing
multiplicities) of terminal quotient singularities of type
$\frac{1}{r_i}(1,-1,b_i)$, $i\in I$ where $I$ is a finite index set.
For simplicity, we will always denote a terminal quotient
singularity $\frac{1}{r}(1,-1,b)$ by  a pair of integers $(b,r)$.
So we will write a basket as:
$$B:=\{n_i\times (b_i,r_i)|i\in J,\ n_i\in {\mathbb Z}^+\},$$
where $n_i$ denotes the multiplicities.

Given baskets $B_1=\{n_i\times (b_i,r_i)\}$ and $B_2=\{m_i\times
(b_i,r_i)\}$, we define
$$B_1 \cup B_2:=\{ (n_i+m_i) \times (b_i,r_i)\}.$$

\begin{defn}
A generalized basket means a collection of pairs of
integers $(b,r)$ with $0<b<r$, not necessarily coprime and allowing
multiplicities.
\end{defn}

\numberwithin{equation}{section}

\begin{setup}\label{inv_bas}{\bf Invariants of baskets.} Given a generalized  basket
$(b,r)$  with $b\leq \frac{r}{2}$ and a fixed integer $n>0$. Let
$\delta:=\lfloor \frac{bn}{r} \rfloor$. Then $\frac{\delta+1}{n}
> \frac{b}{r} \ge \frac{\delta}{n} $.  We define

\begin{equation} \Delta^n({b,r}):=
\delta bn-\frac{(\delta^2+\delta)}{2}r.
\end{equation}

One can see that $\Delta^n({b,r})$ is a non-negative integer. For
a generalized basket $B= \{(b_i,r_i)\}_{ i\in I}$ and a fixed $n>0$,
we define $\Delta^n(B):= \underset{i\in I}\sum
\Delta^n({b_i,r_i})$. By definition, $\Delta^2(B)=0$ for any
basket $B$. By a direct calculation, one gets the following
relation:
$$\frac{\overline{jb_i}(r_i-\overline{jb_i})}{2r_i} - \frac{{jb_i}(r_i-{jb_i})}{2r_i} =
\Delta^j({b_i,r_i})$$ for all $j>0$, $i\in I$. Define

\begin{equation} \sigma(B): = \underset{i\in I}\sum b_i\ \
\text{and}\ \ \sigma'(B):=\underset{i\in I}\sum \frac{b_i^2}{r_i}.
\end{equation}
\end{setup}

\begin{setup}\label{pk}{\bf Packing.} Given a generalized basket
$$B=\{(b_1,r_1),(b_2,r_2),\cdots, (b_k,r_k)\},$$  we call the basket
$$B':=\{ (b_1+b_2,r_1+r_2), (b_3,r_3), \cdots, (b_k,r_k)
\}$$ a packing of $B$ (and $B$ is an unpacking of $B'$), written
as $B \succ B'$. (The symbol $B\succcurlyeq B'$ means either
$B\succ B'$ or $B=B'$.)

If,
furthermore, $b_1r_2-b_2r_1 =1$, we call $B \succ B'$ {\it a prime
packing}. A prime packing is said to have {\it level $n$} if $r_1+r_2=n$.
\end{setup}

The seemingly mysterious notion of packings can  indeed be
realized in various elementary birational maps.

\begin{exmp} We consider the Kawamata blowup \cite{Ka2}.
Let $X=X_\Sigma$ be a toric threefold associated to the fan
$\Sigma$. Suppose that there is a cone $\sigma$ in $\Sigma$
generated by $v_1=(1,0,0),v_2=(0,1,0)$ and $v_3=(s,r-s,1)$ with $
0 <s <r$ and $(s,r)=1$. The cone $\sigma$ gives rise to a quotient
singularity $P \in X$ of type $\frac{1}{r}(r-s,s,1)$.

Let $\pi:\tilde{X} \to X$ be the partial resolution obtained by
the subdivision by adding $v_4=(1,1,1)$. One sees that $\tilde{X}$
has two quotient singularities of type
$\frac{1}{s}(\overline{r},\overline{-r},1)$, and
$\frac{1}{r-s}(\overline{r},\overline{-r},1)$ respectively, where
$\bar{.}$ denotes the residue modulo $s$ and $r-s$ respectively.

 Then it's easy to verify that
$B(X)=\{(b,r)\}$ and $B(\tilde{X})=\{(b',s),(b-b',r-s)\}$ for some
$b,b'$ satisfying  $b'r-bs= \pm 1$. One sees that
$$B(\tilde{X}) \succ B(X)$$ is a prime packing of baskets.
\end{exmp}

\begin{exmp}
Let $X=X_\Sigma$ be a toric threefold associated to the fan
$\Sigma$. Suppose that  there are  two cones $\sigma_4,\sigma_3$
in the fan $\Sigma$ such that
$$\begin{array}{l} \sigma_4 \text{ is generated by } v_1=(1,0,0),v_2=(0,1,0), v_3=(0,0,1) \\
 \sigma_3 \text{ is generated by } v_1=(1,0,0),v_2=(0,1,0), v_4=(s,r-s,-1). \end{array}$$
with $ 0 <s <r$ and $(s,r)=1$.

Let $X^+$ be the threefold obtained by replacing $\sigma_4,
\sigma_3$ with $\sigma_1,\sigma_2$ that
$$\begin{array}{l} \sigma_1 \text{ is generated by } v_2,v_3,v_4 \\
 \sigma_2 \text{ is generated by } v_1,v_3,v_4. \end{array}$$

The birational map $X \dashrightarrow X^+$ is a toric flip. One
can verify that $B(X)=\{(b,r)\}$ and
$B({X^+})=\{(b',s),(b-b',r-s)\}$ for some $b,b'$ satisfying
$b'r-bs= \pm 1$. Similarly,
$$B({X^+}) \succ B(X)$$
is again a prime packing of baskets.

\end{exmp}

We have the following basic properties.
\begin{lem} \label{packing} Let $B \succ B'$ be any packing between generalized baskets.
Keep the same notation as above. Then:
\begin{itemize}
\item[(1)]
$\Delta^n(B) \ge \Delta^n(B')$ for all $n \ge 2$;
\item[(2)] the equality in (1) holds if and only if both $
\frac{b_1}{r_1}$ and $\frac{b_2}{r_2}$ are in the closed interval $[\frac{\delta}{n}, \frac{\delta+1}{n}]$ for some
$\delta$;
\item[(3)]$\sigma(B')=\sigma(B)$ and $\sigma'(B) = \sigma'(B') +
\frac{(r_1b_2-r_2b_1)^2}{r_1r_2(r_1+r_2)} \geq
\sigma'(B').$ Thus equality holds only when
$\frac{b_1}{r_1}= \frac{b_2}{r_2}$.
\end{itemize}
\end{lem}

\begin{proof} First, if both $
\frac{b_1}{r_1}$ and $\frac{b_2}{r_2}$ are in the closed interval $[\frac{\delta}{n}, \frac{\delta+1}{n}]$ for some
$\delta$, then a direct calculation
shows $\Delta^n(B) = \Delta^n(B')$.

Suppose, for some $\delta>j$,
$$\frac{\delta+1}{n}> \frac{b_2}{r_2} \geq \frac{\delta}{n} \geq
\frac{j+1}{n} > \frac{b_1}{r_1} \geq \frac{j}{n}$$ and
 $\frac{j_1+1}{n}>\frac{b_1+b_2}{r_1+r_2} \geq \frac{j_1}{n}$
for some $j_1\in [j,\delta]$. Then
\begin{eqnarray*}
\Delta^n({b_1+b_2,r_1+r_2})&=&
j_1n(b_1+b_2)-\frac{1}{2}(j_1^2+j_1)(r_1+r_2)\\
&=& \Delta^n({b_2,r_2})+\Delta^n({b_1,r_1})+\nabla_2+\nabla_1,
\end{eqnarray*}
where
$\nabla_2=(j_1-\delta)nb_2+\frac{1}{2}(\delta^2+\delta-j_1^2-j_1)r_2$
and $\nabla_1=(j_1-j)nb_1+\frac{1}{2}(j^2+j-j_1^2-j_1)r_1$. Now
since $nb_2\geq \delta r_2$, one gets
$$\nabla_2\leq
\frac{1}{2}(\delta-j_1)(j_1+1-\delta)r_2.$$ When $j_1=\delta$,
$\nabla_2=0$; when $j_1=\delta-1$, $\nabla_2=-nb_1+\delta r_2\leq
0$; when $j_1<\delta-1$, $\nabla_2<0$.

Similarly the relation $nb_1<(j+1)r_1$ implies
$$\nabla_1\leq \frac{1}{2}(j_1-j)(j+1-j_1)r_1.$$
When $j_1=j$, $\nabla_1=0$; when $j_1=j+1$,
$\nabla_1=nb_1-(j+1)r_1<0$; when $j_1>j+1$, $\nabla_1<0$.

Thus in any case, we see $\Delta^n(B) \geq \Delta^n(B')$, which
implies (1). Furthermore we see $\Delta^n(B) = \Delta^n(B')$ if,
and  only if, $\nabla_2=\nabla_1$=0; if, and only if, $j_1=j$ and
$\delta=j_1+1=j+1$. We have proved (2).

The inequality (3) is obtained by a direct calculation.
\end{proof}

\begin{cor}\label{addition} If $B=\{m\times (b,r)|\ b\leq \frac{r}{2},\ b \text{ coprime to } r\}$
and $B'=\{(mb, mr)\}$ for an integer $m>1$, then
\begin{itemize}
\item[(i)] $\sigma(B')=\sigma(B)$; $\sigma'(B')=\sigma'(B)$;
\item[(ii)] $\Delta^n(B')=\Delta^n(B)$ for any $n>0$.
\end{itemize}
\end{cor}
\begin{proof} This can be obtained by the definition of $\sigma$
and Lemma \ref{packing}.
\end{proof}

\begin{rem} The additive properties in Corollary \ref{addition}
allow us to regard the generalized single basket $\{(mb,mr)\}$ as a
basket $\{m\times (b,r)\}$.
\end{rem}

Besides, a prime packing has the following property:

\begin{lem}\label{combini}
Let $B=\{(b_1,r_1), (b_2, r_2)\} \succ \{(b_1+b_2, r_1+r_2)\}=B'$ be
a prime packing as in \ref{pk}, i.e. $b_1r_2-b_2 r_1 =1$.
Then $$\Delta^{r_1+r_2}({b_1+b_2, r_1+r_2})=\Delta^{r_1+r_2}({b_1,
r_1})+\Delta^{r_1+r_2}({b_2, r_2})-1.$$
\end{lem}

\begin{proof}
 When $b_1r_2-b_2r_1=1$, since $r_1> 1, r_2>  1$, one
 has
$$\frac{b_1+ b_2+1}{ r_1+ r_2} >
\frac{b_1}{r_1}
> \frac{ b_1+ b_2}{ r_1+ r_2} > \frac{b_2}{r_2}>
 \frac{ b_1+ b_2-1}{ r_1+ r_2}.$$
  We set $n=r_1+ r_2$. A direct calculation gives the equality
$$\Delta^{n}({b_1+b_2, r_1+r_2})=\Delta^{n}({b_1,
r_1})+\Delta^{n}({b_2, r_2})-1.$$
\end{proof}

\begin{setup}\label{lim}{\bf Initial basket and limiting process.}
Given a basket $B=\{(b_j,r_j)| b_j \text{ coprime
to } r_j, b_j\leq \frac{r_j}{2}\}_{j \in J}$, we
define a sequence of baskets $\{\mathscr{B}^{(n)}(B)\}$ as
follows.

Take the set $S^{(0)}:=\{\frac{1}{n}\}_{n \ge 2}$. For any element
$B_j=(b_j,r_j) \in B$, we can find a unique $n>0$ such that
$\frac{1}{n} > \frac{b_j}{r_j} \geq \frac{1}{n+1}$.
The element $(b_j,r_j)$ can be regarded as finite step successive
packings beginning from the basket $B_j^{(0)}:=\{(nb_j+b_j-r_j)
\times (1,n), (r_j-nb_j) \times (1,n+1)\}$. Adding up those
$B_j^{(0)}$, one obtains the basket $\mathscr{B}^{(0)}(B)=\{
n_{1,2} \times (1,2), n_{1,3} \times (1,3),\cdots, n_{1,r} \times
(1,r)\}$, called {\it the initial basket} of $B$. Clearly
$\mathscr{B}^{(0)}(B)\succcurlyeq B$. Defined in this way,
$\mathscr{B}^{(0)}(B)$ is uniquely determined by the given basket
$B$.

We begin to construct other baskets $\{\mathscr{B}^{(n)}(B)\}$ for
$n> 1$. Consider the sets $S^{(4)}=S^{(3)}=S^{(2)}=S^{(1)}=S^{(0)}$
and
$$S^{(5)}:=S^{(0)} \cup \{\frac{2}{5}\}$$
and inductively, $S^{(n)}=S^{({n-1})} \cup \{\frac{i}{n}
\}_{i=2,\cdots,\lfloor \frac{n}{2} \rfloor}$. Reordering elements in
$S^{(n)}$ and writing $S^{(n)}=\{w^{(n)}_i\}_{i\in I}$
 such that $w^{(n)}_i > w^{(n)}_{i+1}$ for all $i$, then we see that
 the interval $(0,\frac{1}{2}]=\cup_i[w^{(n)}_{i+1}, w^{(n)}_i]$.
 Note that $w^{(n)}_i=\frac{q_i}{p_i}$ with  $p_i$ coprime to $q_i$ and
 $p_i \leq n$ unless $w^{(n)}_i=\frac{1}{m}$ for some $m >n$. First
 we prove the following:
\medskip

\noindent {\bf Claim A.} {\em $u_{1}v_2-u_2v_{1}=1$ for any two
endpoints of
 $[w^{(n)}_{i+1}, w^{(n)}_i]=[\frac{v_{1}}{u_{1}}, \frac{v_2}{u_2}]$.}

\begin{proof}
We can prove this inductively. Suppose that this property holds for
$S^{({n-1})}$. Now, for any $\frac{j}{n} \in S^{(n)}-S^{(n-1)}$,
$\frac{j}{n} \in [w^{(n-1)}_{i+1},w^{(n-1)}_i]=[ \frac{q_{1}}{p_{1}}, \frac{q_{2}}{p_{2}}]$ for some $i$. Thus $
\frac{q_{1}}{p_{1}} < \frac{j}{n} < \frac{q_{2}}{p_{2}}$. If
$p_2 \ge n$, then $\frac{q_2}{p_2}=\frac{1}{m}$ and
$\frac{q_{1}}{p_{1}}=\frac{1}{m+1}$ for some $m \ge n$ which
contradicts to $\frac{j}{n} < \frac{q_2}{p_2}$. Therefore, we must
have  $p_{2} <n$. Then we consider $\frac{j-q_2}{n-p_2}$ and it's
easy to see that
$$
\frac{q_{1}}{p_{1}}\le \frac{j-q_2}{n-p_2} < \frac{j}{n} <
\frac{q_{2}}{p_{2}}.$$ Clearly, $\frac{j-q_2}{n-p_2} \in S^{(n-1)}$
and hence $\frac{j-q_2}{n-p_2}=\frac{q_{1}}{p_{1}}$. It follows
that $n=p_2+\alpha p_{1}, j=q_2+\alpha q_{1}$ for some integer
$\alpha>0$.

If $\alpha \ge 2$, then $\frac{q_{1}}{p_{1}}<
\frac{q_2+(\alpha-1)q_{1}}{p_2+(\alpha-1)p_{1}}  <
\frac{q_{2}}{p_{2}}$, and
$\frac{q_2+(\alpha-1)q_{1}}{p_2+(\alpha-1)p_{1}}  \in
S^{(n-1)}$, which is absurd. Thus $\alpha=1$ and then $n=p_2+
p_{1}, j=q_2+ q_{1}$. It's then clear  that $\frac{j}{n}$
is the only element of $S^{(n)}$ inside the interval
$[\frac{q_{1}}{p_{1}},\frac{q_{2}}{p_{2}}]$. Moreover, $j
p_{1}-n q_{1}=1, n q_2 -j p_2=1$. This completes the proof of
the claim.
\end{proof}

Now for a element $B_i=(b_i,r_i) \in B$, if $\frac{b_i}{r_i} \in
S^{(n)}$, then we set $B^{(n)}_i:=\{(b_i,r_i)\}$. If
$\frac{b_i}{r_i}\not\in S^{(n)}$, then $\frac{q_1}{p_1} <
\frac{b_i}{r_i} < \frac{q_2}{p_2}$ for some interval
$[\frac{q_1}{p_1},\frac{q_2}{p_2}]$ due to $S^{(n)}$.  In this
situation, we can unpack $(b_i,r_i)$ to $B^{(n)}_i:=\{(r_i
q_2-b_ip_2) \times (q_1,p_1),(-r_i q_1+b_i p_1) \times
(q_2,p_2)\}$. Adding up those $B^{(n)}_i$, we get a new basket
$\mathscr{B}^{(n)}(B)$. Clearly $\mathscr{B}^{(n)}(B)$ is uniquely
determined according to our construction and $\mathscr{B}^{(n)}(B)
\succcurlyeq B$ for all $n$.
\end{setup}

\noindent{\bf Claim B.} {\em
$\mathscr{B}^{(n-1)}(B)=\mathscr{B}^{(n-1)}(\mathscr{B}^{(n)}(B))
\succcurlyeq \mathscr{B}^{(n)}(B)$ for all $n\geq 1$.}
\begin{proof} Since we have already seen $\mathscr{B}^{(n-1)}(\mathscr{B}^{(n)}(B))
\succcurlyeq \mathscr{B}^{(n)}(B)$ by definition, it suffices to show the
first equality of the claim. By the definition of
$\mathscr{B}^{(n)}$, we only need to verify the statement for each
element $B_i=\{(b_i,r_i)\} \subset B$ and for $n \ge 5$.

If $\frac{b_i}{ r_i} \in S^{(n-1)} \subset S^{(n)}$, then there is
nothing to prove since the equality follows from the definition of
$\mathscr{B}^{(n)}$ and $\mathscr{B}^{(n-1)}$.

If $\frac{b_i}{ r_i} \in  S^{(n)}-S^{(n-1)}$, then this is also
clear since $\mathscr{B}^{(n)}(B_i) = B_i$.

Suppose finally that  $\frac{b_i}{ r_i} \not \in  S^{(n)}$. Then
$\frac{q_1}{p_1} < \frac{b_i}{r_i} < \frac{q_2}{p_2}$ for some
$\frac{q_1}{p_1}=w^{(n)}_{i+1}$ and $\frac{q_2}{p_2}=w^{(n)}_i$.

{\bf Subcase (i).}  If both of $\frac{q_1}{p_1},\ \frac{q_2}{p_2}$
are in $S^{(n)}- S^{(n-1)}$, then $p_1=p_2=n$ and hence
$p_1q_2-p_2q_1 \neq 1$, a contradiction to Claim A.

{\bf Subcase (ii).} If both $\frac{q_1}{p_1}$ and $\frac{q_2}{p_2}$
are in $S^{(n-1)}$, then by definition
$$\mathscr{B}^{(n-1)}(B_i)=\mathscr{B}^{(n)}(B_i)=\mathscr{B}^{(n-1)}(\mathscr{B}^{(n)}(B_i)).$$

{\bf Subcase (iii).} We are left to consider the situation that one
of $\frac{q_1}{p_1},\ \frac{q_2}{p_2}$ is in $S^{(n-1)}$, but
another one is in $S^{(n)}- S^{(n-1)}$. Let us assume, for example,
$\frac{q_1}{p_1}=w^{(n-1)}_{j+1} \in S^{(n-1)}$. Then
$\frac{q_2}{p_2} < w^{(n-1)}_j = \frac{q}{p}\in S^{(n-1)}$. The
proof for the other case is similar. Notice that by the proof of
Claim A, we have $q_2=q_1+q, p_2=p_1+p$. By definition,
\begin{eqnarray*}
&&\mathscr{B}^{(n)}(B_i)=\{(r_i q_2-b_ip_2) \times (q_1,p_1),(-r_i
q_1+b_i
p_1) \times (q_2,p_2)\},\\
&&\mathscr{B}^{(n-1)}(B_i)=\{(r_i q-b_ip) \times (q_1,p_1),(-r_i
q_1+b_i p_1) \times (q,p)\}.
\end{eqnarray*}
Since $\mathscr{B}^{(n-1)}(q_2,p_2)=\{(q_1,p_1),(q,p)\}$, we get the
following by the direct computation: {\small
\begin{eqnarray*}
\mathscr{B}^{(n-1)}(\mathscr{B}^{(n)}(B_i))&=& \{(r_i q_2-b_ip_2)
\times
(q_1,p_1)\}\cup\{(-r_i q_1+b_ip_1) \times (q_1,p_1),\\
&&(-r_i q_1+b_ip_1) \times (q,p)\}\\
&=&\{ (r_i q-b_ip) \times (q_1,p_1),(-r_i q_1+b_i p_1) \times
(q,p)\}.
\end{eqnarray*}}
So we can see
$\mathscr{B}^{(n-1)}(B_i)=\mathscr{B}^{(n-1)}(\mathscr{B}^{(n)}(B_i))$.
We are done.
\end{proof}

By Claim B, we have a sequence $\{\mathscr{B}^{(n)}(B)\}$ of baskets
with the following relation: {\small \begin{equation}
\mathscr{B}^{(0)}(B)=\ldots=\mathscr{B}^{(4)}(B) \succcurlyeq
\mathscr{B}^{(5)}(B) \succcurlyeq \cdots \succcurlyeq \mathscr{B}^{(n)}(B)
\succcurlyeq \cdots \succcurlyeq B.
\end{equation}}
 Clearly, by definition,
$B=\mathscr{B}^{(w)}(B)$ for some $w \gg 0$ for a given finite
basket $B$. Thus, in some sense, $B$ can be realized as the limit
of the sequence $\{\mathscr{B}^{(n)}(B)\}$, which is called {\it
the canonical sequence of $B$}.

Another direct consequence of Claim B is the following property:
\begin{equation}
\mathscr{B}^{(i)}(\mathscr{B}^{(j)}(B))=\mathscr{B}^{(i)}(B)
\end{equation} for $i\leq j$.

\begin{setup}{\bf The quantity $\epsilon_n(B)$.} Now let us consider the step
$\mathscr{B}^{(n-1)}(B) \succ \mathscr{B}^{(n)}(B)$. For an element
$w \in S^{(n)}$, let $m(w)$ be the number of basket $(b,r)$ in
$\mathscr{B}^{(n)}(B)$ with $b$ coprime to $r$ and $\frac{b}{r}=w$.
Thus we can write $\mathscr{B}^{(n)}(B)=\{m(w) \times
(b,r)\}_{w=\frac{b}{r} \in S^{(n)}}$.

Suppose that $S^{(n)}-S^{(n-1)}=\{\frac{j_s}{n}\}_{s=1,\cdots,t}$.
We have $w^{(n-1)}_{i_s}=\frac{q_{i_s}}{p_{i_s}} > \frac{j_s}{n}
> w^{(n-1)}_{i_s+1}=\frac{q_{i_s+1}}{p_{i_s+1}} $ for some $i_s$.  We remark
that by the proof of Claim A, $j_s= q_{i_s}+q_{i_s+1}$, $n=
p_{i_s}+p_{i_s+1}.$ Since $\mathscr{B}^{(n-1)}(B) =
\mathscr{B}^{(n-1)}(\mathscr{B}^{(n)}(B))$ by Claim B, we may write
$$\mathscr{B}^{(n)}(B)=\{m(w) \times (b,r)\}_{w=\frac{b}{r} \in
S^{(n-1)}}\cup\{m(\frac{j_s}{n}) \times (j_s,n)\}_{\frac{j_s}{n}}
.$$
Then
\begin{eqnarray*}
\mathscr{B}^{(n-1)}(B)=&&\{m(w) \times (b,r)\}_{w=\frac{b}{r} \in
S^{(n-1)}}\cup\{m(\frac{j_s}{n}) \times (q_{i_s},p_{i_s}),\\
&&m(\frac{j_s}{n}) \times (q_{i_s+1},p_{i_s+1})\}_{\frac{j_s}{n}}.
\end{eqnarray*}

We define $\epsilon_n(B):= \sum_{s=1}^t m(\frac{j_s}{n}),$ which is
the number of type $(j_s,n)$ single baskets with $\frac{j_s}{n} \in
S^{(n)}-S^{(n-1)} $. In other words, $\epsilon_n(B)$ counts the
 number of elements $\{(j_s,n)\}$ contained in
$\mathscr{B}^{(n)}(B)$ with $(j_s,n)=1$ and $j_s
>1$. By Claim A, we  conclude that $\mathscr{B}^{(n-1)}(B) \succcurlyeq \mathscr{B}^{(n)}(B)$ consists of
$\epsilon_n(B)$ prime packings of level $n$.  This is going to be
an important quantity in our arguments.
\end{setup}

\begin{defn}
Given a basket $B$. The sequence defined as in $(2.3)$ is called the
canonical sequence of prime unpackings of $B$, or canonical sequence
of $B$ for short.
\end{defn}


\begin{setup} {\bf Notation.} When no confusion is likely, we will
simply write $B^{(n)}$ for $\mathscr{B}^{(n)}(B)$.

\end{setup}

\begin{lem}\label{delta} For the canonical sequence $\{B^{(n)}\}$ , the following statements
hold.
\begin{itemize}
\item[(i)]
$\Delta^j(B^{(0)})= \Delta^j(B)$ for $j=3,4$;
\item[(ii)]
$\Delta^j(B^{(n-1)})= \Delta^j(B^{(n)})$ for all $j <n$;
\item[(iii)]
$\Delta^n(B^{(n-1)})= \Delta^n(B^{(n)})+\epsilon_n(B)$.
\item[(iv)]
$\Delta^n(B^{(n)})=\Delta^n(B)$.
\end{itemize}
\end{lem}

\begin{proof}
From $B^{(0)}$ to $B$, via $B^{(n)}$, the whole process can be
realized through a composition of finite number of prime
packings. Each step is of the form $\{ (q_1,p_1),(q_2,p_2) \} \succ
\{(q_1+q_2,p_1+p_2)\}$. Notice that either
$\frac{q_1}{p_1},\frac{q_2}{p_2} \leq \frac{1}{3}$ or
$\frac{q_1}{p_1},\frac{q_2}{p_2} \geq \frac{1}{3}$. By Lemma
\ref{packing}(2), one gets $\Delta^3(B^{(0)})=\Delta^3(B)$. The
proof for $\Delta^4$ is similar.

Now we consider the typical step $B^{(n-1)} \succ B^{(n)}$. By Lemma
\ref{combini} and a direct computation, one has:
$$\begin{array} {l}\Delta^n(B^{(n-1)})-\Delta^n( B^{(n)})\\
= \sum_{s=1}^t m(\frac{j_s}{n})
(\Delta^n(q_{i_s},p_{i_s}+\Delta^n(q_{i_s+1},p_{i_s+1})-\Delta^n(j_s,n))\\
=\sum_{s=1}^t m(\frac{j_s}{n})
(\Delta^n(q_{i_s},p_{i_s})+\Delta^n(q_{i_s+1},p_{i_s+1})-\Delta^n(q_{i_s}+q_{i_s+1},p_{i_s}+p_{i_s+1}))\\
 =\sum_{s=1}^t m(\frac{j_s}{n})\\
=\epsilon_n(B), \end{array}$$ where one notices
$n=p_{i_s}+p_{i_s+1}$.

Finally, for any $j <n$, and suppose that $\frac{k+1}{j} \ge
\frac{q_{i_s}}{p_{i_s}}=w^{(n-1)}_{i_s} > \frac{k}{j}$ for some
$k$. Then $\frac{k+1}{j} \in S^{(n-1)}$ by definition. Thus
$\frac{q_{i_s+1}}{p_{i_s+1}}=w^{(n-1)}_{i_s+1} \ge \frac{k}{j}$. By
Lemma \ref{packing}, we have
$$
\Delta^j(q_{i_s},p_{i_s})+\Delta^j(q_{i_s+1},p_{i_s+1})=\Delta^j(q_{i_s}+q_{i_s+1},p_{i_s}+p_{i_s+1}).$$

The last statement is due to (ii) and the fact that $B=B^{(n)}$ for
a sufficiently large $n$. This completes the proof.
\end{proof}

Let us go back to investigate the canonical sequence (2.3)
$$ B^{(0)} \succcurlyeq B^{(5)} \succcurlyeq ... \succcurlyeq B^{(n)} \succcurlyeq ... \succcurlyeq B.$$
We see that $\Delta^j (B^{(n)})= \Delta^j (B)$ for all $j <n$.
Thus we can informally view $B^{(n)}$  as an $n$-th order
approximation of $B$. Also each approximation step $B^{(n-1)}
\succcurlyeq B^{(n)}$ is nothing but the composition of prime packings of
$\epsilon_n$ pairs of baskets of type $(b,n)$ with $b$ coprime to
$n$, $b\leq \frac{r}{2}$ and $b>1$.
\medskip

\section{\bf Formal baskets}
In this section, we are going to introduce the notion of formal
baskets. A formal basket is a basket together with a choice of  $K^3$ and $\chi$.
The purpose of this section is to
classify all formal baskets with a given initial sequence $(\chi_1,\ldots,\chi_k)$.

Given a  3-fold $X$ with canonical singularities, there is an
associated basket $B:={\mathscr{B}}(X)$\footnote{Iano-Fletcher
\cite{Fletcher} has shown that Reid's virtual basket
$\mathscr{B}(X)$ is uniquely determined by $X$.} according to
Reid.

\begin{setup}{\bf Euler characteristic.}\label{pm}
Let us recall Reid's Riemann-Roch formula (\cite[Page 143]{YPG}) for
a $\bQ$-factorial terminal 3-fold $X$:
 for all $m>1$,
$$\chi(X,\mathcal{O}_X(mK_X))=\frac{1}{12}m(m-1)(2m-1)K_X^3-(2m-1)\chi(\mathcal
{O}_X)+l(m)  \eqno (3.1)$$ where the correction term $l(m)$ can be
computed as:
$$l(m):=\sum_{Q\in \mathscr{B}(X)}l_Q(m):=\sum_{Q\in \mathscr{B}(X)}\sum_{j=1}^{m-1}
\frac{\overline{jb_Q}(r_Q-\overline{jb_Q})}{2r_Q}$$ where the sum
$\sum_{Q}$ runs through all single baskets $Q$ in $\mathscr{B}(X)$
with type $\frac{1}{r_{Q}}(1,-1,b_{Q})$ and $\overline{jb_Q}$ means
the smallest residue of $jb_Q$ \text{mod} $r_Q$.

For brevity, $\chi(X,\mathcal{O}_X(mK_X))$ is usually denoted as
$\chi_m(X)$ or simply $\chi_m$.
\medskip

We are going to analyze the above formula and Reid's virtual basket
$\mathscr{B}(X)$.
\end{setup}

\begin{setup}\label{del}{\bf  Euler characteristic in terms of baskets.}
Take $B={\mathscr{B}}(X)$ and set
$\Delta:=\Delta(B),\sigma:=\sigma(B),\sigma':=\sigma'(B)$ (cf.
2.2). We can now rewrite Reid's Riemann-Roch formula as
the following:
$$\left\{ \begin{array}{ll} \chi_2 &= \frac{1}{2} (K_X^3-\sigma') + \frac{1}{2} \sigma-3 \chi,\\
 \chi_3-\chi_2 & = \frac{4}{2}  (K_X^3-\sigma')  + \frac{2}{2} \sigma-2 \chi,\\
 \chi_{m+1}-\chi_m & = \frac{m^2}{2}(K_X^3-\sigma')  + \frac{m}{2} \sigma-2 \chi
+\Delta^m, \text{ for } m \geq 3.
\end{array} \right. \eqno(3.2)
$$
\end{setup}

Notice that, by the equalities $(3.2)$, all $\chi_m$ are
determined by $\sigma, \sigma'-K^3, \chi, \Delta^j$ for all $j
<m$. These, in turn, are determined by $B, \chi$ and $\chi_2$ by
virtue of the first equality in $(3.2)$. This leads us to consider
a more general setting.

\begin{defn} Assume that $B$ is a basket, $\tilde{\chi}$ and  $\tilde{\chi_2}$  are integers. We call the triple  ${\bf B}:=(B, \tilde{\chi},
\tilde{\chi_2})$ a {\it formal basket}.
\end{defn}

We can define the Euler characteristic and $K^3$ of a formal
basket formally by the Riemann-Roch formula. First we
define
$$\left\{ \begin{array}{l} \chi_2({\bf B}):=\tilde{\chi}_2,\\ \chi_3({\bf B}):=-\sigma(B)+ 10
\tilde{\chi}+5\tilde{\chi}_2 \end{array} \right.$$ and the volume
$$\begin{array}{lll} K^3({\bf B})&:=&\sigma'({B})-4 \tilde{\chi}
-3\tilde{\chi}_2+\chi_3({\bf
B})\\
&=& -\sigma+\sigma'+6\tilde{\chi}+2\tilde{\chi}_2.
\end{array} \eqno(3.3)
$$ For $m \geq 4$, the Euler characteristic
$\chi_m({\bf B})$  is defined inductively by {\small
$$\chi_{m+1}({\bf B})-\chi_m({\bf B}):=
\frac{m^2}{2}(K^3({\bf B})-\sigma'({B})) +\frac{m}{2}
\sigma({B})-2\tilde{\chi}+\Delta^m(B).\eqno{(3.4)}$$}

Clearly, by definition, $\chi_m({\bf B})$ is an integer for all
$m\geq 4$ because $K^3({\bf B})-\sigma'({B})=-4 \tilde{\chi}
-3\tilde{\chi}_2+\chi_3({\bf B})$ and $\sigma=10
\tilde{\chi}+5\tilde{\chi}_2-\chi_3({\bf B})$ have the same
parity.


Given a $\mathbb{Q}$-factorial canonical 3-fold $X$, one can associate to $X$ a triple ${\bf
B}(X):=(B,\tilde{\chi}, \tilde{\chi_2}) $ where $B=\mathscr{B}(X)$,
$\tilde{\chi}=\chi(\mathcal{O}_X)$ and $\tilde{\chi_2}= \chi_2(X)$. It's
clear that such a triple is a formal basket. The Euler characteristic  and $K^3$ of the formal basket ${\bf B}(X)$ are nothing but
the Euler characteristic  and $K^3$ of the variety $X$.

\begin{setup}{\bf Notations.}
For simplicity, we denote $\chi_{m}({\bf B})$ by $\tilde{\chi}_m$
for all $m\geq 2$. Also denote $K^3({\bf B})$ by $\tilde{K}^3$,
$\sigma=\sigma(B)$, $\sigma'=\sigma'(B)$ and
$\Delta^m=\Delta^m(B)$.
\end{setup}

\begin{defn} Let ${\bf B}:=(B, \tilde{\chi},
\tilde{\chi_2})$ and ${\bf B'}:=(B', \tilde{\chi},
\tilde{\chi_2})$ be two formal baskets.

(1) We say that ${\bf B'}$ is a packing of ${\bf B}$ (written as
${\bf B}\succ {\bf B'}$) if $B\succ B'$. Clearly ``packing'' between
formal baskets gives a partial ordering.

(2) A formal basket ${\bf B}$ is called {\it positive} if $K^3({\bf
B})>0$.

(3) A formal basket ${\bf B}$ is said to be minimal positive if it
is positive and minimal with regard to packing relation.
\end{defn}

By definition and Lemma \ref{packing}(1), we immediately get the
following:

\begin{lem}\label{inequal} Assume ${\bf B}:=(B, \tilde{\chi},
\tilde{\chi_2})\succ {\bf B'}:=(B', \tilde{\chi},
\tilde{\chi_2})$. Then
\begin{itemize}
\item[(1)] $K^3({\bf B})\geq K^3({\bf B'})$;

\item[(2)] $\chi_m({\bf B})\geq \chi_m({\bf B'})$ for all $m\geq 2$.
\end{itemize}
\end{lem}

In what follows, we would like to classify
all baskets with a given initial sequence
$(\tilde{\chi},\tilde{\chi_2},\tilde{\chi_3},\cdots,\tilde{\chi_m})$.

First of all, by the definition of $\tilde{K}^3$ and
$\tilde{\chi}_m$, we get:
$$\begin{array}{rl} \tau:=\sigma'-\tilde{K}^3&=   4 \tilde{\chi}+3\tilde{\chi}_2-\tilde{\chi}_3,\\
\sigma &=  10 \tilde{\chi} +5 \tilde{\chi}_2-\tilde{\chi}_3\\
\Delta^3&=  5 \tilde{\chi} +6\tilde{\chi}_2-4\tilde{\chi}_3+\tilde{\chi}_4\\
\Delta^4&=  14 \tilde{\chi} +14\tilde{\chi}_2-6\tilde{\chi}_3-\tilde{\chi}_4+\tilde{\chi}_5\\
\Delta^5 &=27 \tilde{\chi} +25\tilde{\chi}_2- 10 \tilde{\chi}_3 - \tilde{\chi}_5 + \tilde{\chi}_6  \\
\Delta^6 &= 44 \tilde{\chi} +39\tilde{\chi}_2- 15 \tilde{\chi}_3- \tilde{\chi}_6+\tilde{\chi}_7    \\
\Delta^7&= 65 \tilde{\chi} +56\tilde{\chi}_2- 21 \tilde{\chi}_3 - \tilde{\chi}_7 + \tilde{\chi}_8 \\
\Delta^8&= 90 \tilde{\chi} +76\tilde{\chi}_2- 28 \tilde{\chi}_3 - \tilde{\chi}_8 + \tilde{\chi}_9 \\
\Delta^9&= 119 \tilde{\chi} +99\tilde{\chi}_2- 36 \tilde{\chi}_3  - \tilde{\chi}_9 +\tilde{\chi}_{10} \\
\Delta^{10}&= 152 \tilde{\chi} +125\tilde{\chi}_2- 45 \tilde{\chi}_3 - \tilde{\chi}_{10} + \tilde{\chi}_{11} \\
\Delta^{11}&= 189 \tilde{\chi} +154\tilde{\chi}_2- 55 \tilde{\chi}_3 - \tilde{\chi}_{11} +\tilde{\chi}_{12}  \\
\Delta^{12}&=  230 \tilde{\chi} +186\tilde{\chi}_2- 66 \tilde{\chi}_3- \tilde{\chi}_{12} +\tilde{\chi}_{13}  \\
\end{array} \eqno(3.5)
$$
Recall that $B^{(0)}=\{n^0_{1,2} \times (1,2),\cdots,
n^0_{1,r} \times (1,r)\}$ is the initial basket of $B$.
Then by Lemma \ref{delta} and the definition of
$\sigma(B)$, we have
$$ \begin{array}{l} \sigma(B)= \sigma(B^{(0)})=\sum n^0_{1,r},\\
\Delta^3(B)= \Delta^3(B^{(0)})= n^0_{1,2} \\
\Delta^4(B)=\Delta^4(B^{(0)})=2 n^0_{1,2}+ n^0_{1,3}
\end{array}$$
Therefore, the initial basket has the coefficients:
$$B^{(0)} \left\{ \begin{array}{l}
n^0_{1,2}=5 \tilde{\chi}+ 6 \tilde{\chi}_2 -4 \tilde{\chi}_3 +\tilde{\chi}_4\\
n^0_{1,3}=4 \tilde{\chi} +2 \tilde{\chi}_2+2\tilde{\chi}_3-3 \tilde{\chi}_4 +\tilde{\chi}_5\\
n^0_{1,4}=\tilde{\chi} -3 \tilde{\chi}_2+\tilde{\chi}_3+2\tilde{\chi}_4-\tilde{\chi}_5- \sum_{r \ge 5} n^0_{1,r}\\
n^0_{1,r}, r \ge 5.
\end{array} \right. \eqno(3.6)
$$
By Lemma \ref{delta}, we see
$$\begin{array}{lll} \epsilon_5&:=&\Delta^5(B^{(0)})-\Delta^{5}(B)=4
n^0_{1,2}+2n^0_{1,3}+n^0_{1,4}-\Delta^{5}(B)\\
&=&2 \tilde{\chi}- \tilde{\chi}_3 +2 \tilde{\chi}_5
-\tilde{\chi}_6-\sigma_5\ \ \text{where} \\
\sigma_5 &:=& \sum_{r \ge 5} n^0_{1,r}.
\end{array} \eqno(3.7)$$
Thus we can write {\footnotesize
$$B^{(5)}=\{ n^5_{1,2} \times (1,2),n^5_{2,5} \times (2,5), n^5_{1,3}
\times (1,3), n^5_{1,4}  \times (1,4), n^5_{1,5} \times
(1,5),\cdots\}$$} with
$$B^{(5)}\left\{
 \begin{array}{l}
  n^5_{1,2}=3 \tilde{\chi} +6\tilde{\chi}_2- 3 \tilde{\chi}_3 + \tilde{\chi}_4 - 2 \tilde{\chi}_5 + \tilde{\chi}_6+ \sigma_5 ,\\
 n^5_{2,5}= 2 \tilde{\chi}-\tilde{\chi}_3 +   2 \tilde{\chi}_5  - \tilde{\chi}_6- \sigma_5 \\
 n^5_{1,3}= 2 \tilde{\chi}+2\tilde{\chi}_2+ 3 \tilde{\chi}_3- 3 \tilde{\chi}_4 -\tilde{\chi}_5 + \tilde{\chi}_6   + \sigma_5 , \\
 n^5_{1,4}= \tilde{\chi} -3\tilde{\chi}_2 +  \tilde{\chi}_3 +2 \tilde{\chi}_4 -\tilde{\chi}_5 -\sigma_5 \\
  n^5_{1,r}=n^0_{1,r}, r \ge 5
\end{array} \right. \eqno(3.8)$$ Noting that this is obtained from $B^{(0)}$ by
taking $\epsilon_5$ prime packings of type  $\{(1,2),
(1,3)\}\succ \{(2,5)\}$.

Clearly, $B^{(5)}=B^{(6)}$ by our construction. Thus by Lemma
\ref{delta} we have
$\Delta^6(B^{(5)})=\Delta^6(B^{(6)})=\Delta^6(B)$. Computation
shows that
$$
\begin{array}{lll}
\Delta^6(B^{(5)})&=&6n^5_{1,2}+9n^5_{2,5}+3n^5_{1,3}+2n^5_{1,4}+n^5_{1,5}\\
&=&
44\tilde{\chi}+36\tilde{\chi}_2-16\tilde{\chi}_3+\tilde{\chi}_4+
\tilde{\chi}_5-\epsilon,
\end{array}$$
where $$\epsilon:=n^0_{1,5}+2 \sum_{r \ge 6} n^0_{1,r} = 2
\sigma_5-n^0_{1,5} \ge 0. \eqno(3.9)$$ Compare this with $(3.5)$,  we see
$$
\epsilon_6=-3\tilde{\chi}_2-\tilde{\chi}_3+\tilde{\chi}_4+\tilde{\chi}_5+\tilde{\chi}_6-\tilde{\chi}_7-\epsilon=0.
\eqno(3.10)$$

Next, by similar computation,  we get
$$\begin{array}{ll}\epsilon_7:&=\Delta^7(B^{(6)})-\Delta^7(B)=\Delta^7(B^{(5)})-\Delta^7(B)\\
&=9n^5_{1,2}+13n^5_{2,5}+5n^5_{1,3}+3n^5_{1,4}+2n^5_{1,5}+n^5_{1,6}-\Delta^7(B)\\
&=\tilde{\chi}-\tilde{\chi}_2-\tilde{\chi}_3+\tilde{\chi}_6+\tilde{\chi}_7-\tilde{\chi}_8-2\sigma_5+2n^0_{1,5}+n^0_{1,6}.
\end{array} \eqno(3.11)$$
Since $S^{(7)}-S^{(6)}=\{\frac{2}{7},\frac{3}{7}\}$, there are two
ways of prime packings into type $(b,7)$ baskets. Let $\eta \geq 0$
be the number of prime packings of type $\{(1,3),(1,4)\} \succ
\{(2,7)\}$. Then $\epsilon_7-\eta \ge 0$ is the number of prime
packings of type $\{(1,2),(2,5)\} \succ \{(3,7)\}$. Thus we can
write $B^{(7)}=\{ n^7_{b,r} \times (b,r)\}_{\frac{b}{r} \in
S^{(7)}}$ with {\small
$$ B^{(7)}\left\{
\begin{array}{l}
n^7_{1,2} =2 \tilde{\chi} +7\tilde{\chi}_2 - 2 \tilde{\chi}_3 +\tilde{\chi}_4  - 2 \tilde{\chi}_5 - \tilde{\chi}_7+
\tilde{\chi}_8+3 \sigma_5-2n^0_{1,5}-n^0_{1,6} + \eta  \\
 n^7_{3,7} =\tilde{\chi}-\tilde{\chi}_2-\tilde{\chi}_3+\tilde{\chi}_6+\tilde{\chi}_7-\tilde{\chi}_8-2\sigma_5+2n^0_{1,5}+n^0_{1,6}- \eta\\
 n^7_{2,5} = \tilde{\chi} +\tilde{\chi}_2+ 2 \tilde{\chi}_5 - 2 \tilde{\chi}_6 - \tilde{\chi}_7+ \tilde{\chi}_8 + \sigma_5-2n^0_{1,5}-n^0_{1,6}  + \eta \\
 n^7_{1,3} =2 \tilde{\chi} +2\tilde{\chi}_2+ 3 \tilde{\chi}_3 - 3 \tilde{\chi}_4 -\tilde{\chi}_5 + \tilde{\chi}_6  + \sigma_5  - \eta\\
n^7_{2,7} = \eta\\
n^7_{1,4} = \tilde{\chi} -3\tilde{\chi}_2+ \tilde{\chi}_3 + 2 \tilde{\chi}_4 - \tilde{\chi}_5 - \sigma_5 - \eta\\
n^7_{1,r}=n^0_{1,r}, r \ge 5
\end{array} \right. \eqno(3.12)$$}

From $B^{(7)}$, we can compute $\epsilon_8$ and then $B^{(8)}$, and
inductively $B^{(n)}$ for all $n\geq 9$. But notice that one can
even compute $\epsilon_9$, $\epsilon_{10}$ and $\epsilon_{12}$
directly from $B^{(7)}$, thanks to Lemma \ref{packing}.

To see this, let's consider
$\epsilon_9:=\Delta^9(B^{(8)})-\Delta^9(B)$ for example. Note that
$B^{(7)} \succ B^{(8)}$ is obtained by some prime packings into
$\{(3,8)\}$. Every such packing, which is $\{(2,5),(1,3)\} \succ
\{(3,8)\}$, happens inside a closed interval $[ \frac{3}{9},
\frac{4}{9}]$. Thus by Lemma \ref{packing}(2),
$\Delta^9(B^{(8)})=\Delta^9(B^{(7)})$ and hence
$$ \epsilon_9:=\Delta^9(B^{(8)})-\Delta^9(B)= \Delta^9(B^{(7)})-\Delta^9(B).$$ Similarly we can see
$\Delta^{10}(B^{(9)})=\Delta^{10}(B^{(7)})$ and
$\Delta^{12}(B^{(10)})=\Delta^{12}(B^{(7)})$. Unfortunately,
$\Delta^{11}(B^{(10)})\neq \Delta^{11}(B^{(7)})$.

In summary, we have the following by direct calculations:{\small
$$ \begin{array}{lll}
\Delta^8(B^{(7)})& =&12 n^7_{1,2}+30 n^7_{3,7}+18 n^7_{2,5}+ 7
n^7_{1,3}+11
n^7_{2,7}+4 n^7_{1,4}\\
&&+3 n^7_{1,5}+2 n^7_{1,6}+n^7_{1,7} \\
&=& 90 \tilde{\chi}+74\tilde{\chi}_2 -29\tilde{\chi}_3-\tilde{\chi}_4+\tilde{\chi}_5+\tilde{\chi}_6 -3\sigma_5\\
&&+3 n^0_{1,5}+2 n^0_{1,6}+n^0_{1,7}; \\
\Delta^9(B^{(8)})& =&\Delta^9(B^{(7)})\\
& =&16 n^7_{1,2}+39 n^7_{3,7}+24 n^7_{2,5}+ 9 n^7_{1,3}+15
n^7_{2,7}+6 n^7_{1,4} \\
&&+4 n^7_{1,5}+3n^7_{1,6}+2n^7_{1,7}+n^7_{1,8} \\
&= &119\tilde{\chi}+97\tilde{\chi}_2 -38\tilde{\chi}_3+\tilde{\chi}_4+\tilde{\chi}_5-\tilde{\chi}_7+\tilde{\chi}_8-3 \sigma_5+\eta\\
&&+2 n^0_{1,5}+2n^0_{1,6}+2n^0_{1,7}+n^0_{1,8}; \\
\Delta^{10}(B^{(9)}) &=&\Delta^{10}(B^{(8)})=\Delta^{10}(B^{(7)})\\
& =&20 n^7_{1,2}+50 n^7_{3,7}+30 n^7_{2,5}+ 12 n^7_{1,3}+19
n^7_{2,7}+8 n^7_{1,4}\\
&&+5 n^7_{1,5}+4n^7_{1,6}+3n^7_{1,7}+2n^7_{1,8}+n^7_{1,9} \\
&=& 152\tilde{\chi} +120\tilde{\chi}_2-46\tilde{\chi}_3+2\tilde{\chi}_6-6\sigma_5-\eta \\
&& +5n^0_{1,5}+4n^0_{1,6}+3n^0_{1,7}+2n^0_{1,8}+n^0_{1,9};\\
\Delta^{12}(B^{(11)})&=&\Delta^{12}(B^{(10)})=\cdots=\Delta^{12}(B^{(7)})\\
& =&30 n^7_{1,2}+75 n^7_{3,7}+46 n^7_{2,5}+ 18 n^7_{1,3}+30
n^7_{2,7}+12 n^7_{1,4}\\
&&+9 n^7_{1,5}+6n^7_{1,6}+5n^7_{1,7}+4n^7_{1,8}+3n^7_{1,9}+2n^7_{1,10}+n^7_{1,11} \\
&=& 229\tilde{\chi}+181\tilde{\chi}_2 -69\tilde{\chi}_3+2\tilde{\chi}_5+\tilde{\chi}_6-\tilde{\chi}_7+\tilde{\chi}_8-8 \sigma_5+\eta\\
&&+7 n^0_{1,5}+5n^0_{1,6}+5n^0_{1,7}+4n^0_{1,8}+3n^0_{1,9}+2n^0_{1,10}+n^0_{1,11}. \\
\end{array} $$}

We thus have:{\small
$$
\begin{array}{ll}
 \epsilon_{8} =& -2\tilde{\chi}_2-\tilde{\chi}_3-\tilde{\chi}_4+\tilde{\chi}_5+\tilde{\chi}_6+\tilde{\chi}_8-\tilde{\chi}_9-3\sigma_5\\
&+3 n^0_{1,5}+2 n^0_{1,6}+n^0_{1,7}; \\
\epsilon_{9} =&-2\tilde{\chi}_2-2\tilde{\chi}_3+\tilde{\chi}_4+\tilde{\chi}_5-\tilde{\chi}_7+\tilde{\chi}_8+\tilde{\chi}_9-\tilde{\chi}_{10} -3 \sigma_5+\eta\\
&+2 n^0_{1,5}+2n^0_{1,6}+2n^0_{1,7}+n^0_{1,8}; \\
 \epsilon_{10} =& -5\tilde{\chi}_2-\tilde{\chi}_3+2\tilde{\chi}_6+\tilde{\chi}_{10}-\tilde{\chi}_{11}-6\sigma_5 -\eta\\
 &+5n^0_{1,5}+4n^0_{1,6}+3n^0_{1,7}+2n^0_{1,8}+n^0_{1,9};\\
\epsilon_{12}=
&-\tilde{\chi}-5\tilde{\chi}_2-3\tilde{\chi}_3+2\tilde{\chi}_5+\tilde{\chi}_6-\tilde{\chi}_7+\tilde{\chi}_8+\tilde{\chi}_{12}-\tilde{\chi}_{13}-8
  \sigma_5+\eta\\
  &+7 n^0_{1,5}+5n^0_{1,6}+5n^0_{1,7}+4n^0_{1,8}+3n^0_{1,9}+2n^0_{1,10}+n^0_{1,11}. \\
\end{array} \eqno(3.13)
$$}\\
Since both $\epsilon_{10}$ and $\epsilon_{12}$ are
non-negative, we have $\epsilon_{10}+\epsilon_{12} \geq 0$.
This gives rise to:
$$2\tilde{\chi}_5+3\tilde{\chi}_6+\tilde{\chi}_8+\tilde{\chi}_{10}+\tilde{\chi}_{12} \ge \tilde{\chi} +10
\tilde{\chi}_2+4\tilde{\chi}_3+\tilde{\chi}_7+\tilde{\chi}_{11}+\tilde{\chi}_{13}+R,
\eqno{(3.14)}$$ where {\small
$$R:=14 \sigma_5-12 n^0_{1,5}-9
n^0_{1,6}-8n^0_{1,7}-6n^0_{1,8}-4n^0_{1,9}-2n^0_{1,10}-n^0_{1,11}$$
$$ =
2n^0_{1,5}+5n^0_{1,6}+6n^0_{1,7}+8n^0_{1,8}+10n^0_{1,9}+12n^0_{1,10}+13n^0_{1,11}+14
\sum_{r \ge 12} n^0_{1,r}.$$}

\begin{rem} By definition, $\epsilon_n  \ge 0$. This gives rise to various new inequalities among Euler characteristic.
For example, $\epsilon_5 \ge 0$ (cf. 3.7) gives $$2 \tilde{\chi}- \tilde{\chi}_3 +2 \tilde{\chi}_5
-\tilde{\chi}_6 \ge 0.$$
In particular, for a $\mathbb{Q}$-factorial threefold $X$ with canonical singularities, one has
$ 2 {\chi}(X)- {\chi}_3(X) +2 {\chi}_5(X)
-{\chi}_6(X) \ge 0$.

Among those we have presented above, the equation $(3.10)$ and the
inequality $(3.14)$ will play the most important roles in the
context.
\end{rem}

In practice, we will frequently end up with situations (see Lemma
4.8 and the proof of Theorem 4.12) satisfying the following
assumption and then our computation will be comparatively simpler.

\begin{setup}\label{r6} {\bf Assumption:} $\tilde{\chi}_2=0$ and
$n^0_{1,r}=0$ for all $r \ge 6$.
\end{setup}

Under Assumption \ref{r6}, we list our datum in details as
follows. First,
$$\epsilon_7=\tilde{\chi}-\tilde{\chi}_3+\tilde{\chi}_6+
\tilde{\chi}_7-\tilde{\chi}_8
$$ and
$B^{(7)}=\{ n^7_{b,r} \times (b,r)\}_{\frac{b}{r}
\in S^{(7)}}$
 has coefficients:
$$B^{(7)} \left\{
\begin{array}{l}
n^7_{1,2} =2 \tilde{\chi}  - 2 \tilde{\chi}_3 +\tilde{\chi}_4  - 2 \tilde{\chi}_5 - \tilde{\chi}_7+ \tilde{\chi}_8+ n^0_{1,5} + \eta  \\
 n^7_{3,7} =\tilde{\chi}-\tilde{\chi}_3+\tilde{\chi}_6+\tilde{\chi}_7-\tilde{\chi}_8- \eta\\
 n^7_{2,5} = \tilde{\chi}+ 2 \tilde{\chi}_5 - 2 \tilde{\chi}_6 - \tilde{\chi}_7+ \tilde{\chi}_8 -n^0_{1,5}  + \eta \\
 n^7_{1,3} =2 \tilde{\chi} + 3 \tilde{\chi}_3 - 3 \tilde{\chi}_4 -\tilde{\chi}_5 + \tilde{\chi}_6  + n^0_{1,5}  - \eta\\
n^7_{2,7} = \eta\\
n^7_{1,4} = \tilde{\chi} + \tilde{\chi}_3 + 2 \tilde{\chi}_4 - \tilde{\chi}_5 - n^0_{1,5} - \eta\\
n^7_{1,5}=n^0_{1,5}.
\end{array} \right.$$
We have already known
$$\epsilon_8 =
-\tilde{\chi}_3-\tilde{\chi}_4+\tilde{\chi}_5+\tilde{\chi}_6+\tilde{\chi}_8-\tilde{\chi}_9.$$
Thus, taking some prime packings into account, $B^{(8)}=\{ n^8_{b,r}
\times (b,r)\}_{\frac{b}{r} \in S^{(8)}}$ has the coefficients:
$$B^{(8)} \left\{
\begin{array}{l}
n^8_{1,2} =2 \tilde{\chi}  - 2 \tilde{\chi}_3 +\tilde{\chi}_4  - 2 \tilde{\chi}_5 - \tilde{\chi}_7+ \tilde{\chi}_8+n^0_{1,5} + \eta  \\
 n^8_{3,7} =\tilde{\chi}  -\tilde{\chi}_3 +  \tilde{\chi}_6 + \tilde{\chi}_7- \tilde{\chi}_8  - \eta\\
 n^8_{2,5} = \tilde{\chi} +\tilde{\chi}_3+\tilde{\chi}_4+  \tilde{\chi}_5 - 3 \tilde{\chi}_6 - \tilde{\chi}_7+ \tilde{\chi}_9- n^0_{1,5}  + \eta \\
 n^8_{3,8}=-\tilde{\chi}_3-\tilde{\chi}_4+\tilde{\chi}_5+\tilde{\chi}_6+\tilde{\chi}_8-\tilde{\chi}_9\\
 n^8_{1,3} =2 \tilde{\chi} + 4 \tilde{\chi}_3 - 2 \tilde{\chi}_4 -2\tilde{\chi}_5 - \tilde{\chi}_8 +\tilde{\chi}_9 + n^0_{1,5}  - \eta\\
n^8_{2,7} = \eta\\
n^8_{1,4} = \tilde{\chi} + \tilde{\chi}_3 + 2 \tilde{\chi}_4 - \tilde{\chi}_5 -n^0_{1,5} - \eta\\
n^8_{1,5}=n^0_{1,5}.
\end{array} \right.$$
We know that
$$ \epsilon_9=-2\tilde{\chi}_3+\tilde{\chi}_4+\tilde{\chi}_5
-\tilde{\chi}_7+\tilde{\chi}_8+\tilde{\chi}_9-\tilde{\chi}_{10}
-n^0_{1,5}+\eta.$$ Moreover
$S^{(9)}-S^{(8)}=\{\frac{4}{9},\frac{2}{9}\}$. Let $\zeta$ be the
number of prime packings of type $\{(1,2),(3,7)\} \succ \{(4,9)\}$,
then the number of type $\{(1,4),(1,5)\} \succ \{(2,9)\}$ prime
packings is $\epsilon_9-\zeta$ . We can get $B^{(9)}$ consisting of
the following coefficients.
$$B^{(9)} \left\{
\begin{array}{l}
n^9_{1,2} =2 \tilde{\chi}  - 2 \tilde{\chi}_3 +\tilde{\chi}_4  - 2 \tilde{\chi}_5 - \tilde{\chi}_7+ \tilde{\chi}_8+n^0_{1,5} + \eta-\zeta  \\
n^9_{4,9}=\zeta \\
 n^9_{3,7} =\tilde{\chi}  -\tilde{\chi}_3 +  \tilde{\chi}_6 + \tilde{\chi}_7- \tilde{\chi}_8  - \eta-\zeta\\
 n^9_{2,5} = \tilde{\chi} +\tilde{\chi}_3+\tilde{\chi}_4+  \tilde{\chi}_5 - 3 \tilde{\chi}_6 - \tilde{\chi}_7+ \tilde{\chi}_9- n^0_{1,5}  + \eta \\
 n^9_{3,8}=-\tilde{\chi}_3-\tilde{\chi}_4+\tilde{\chi}_5+\tilde{\chi}_6+\tilde{\chi}_8-\tilde{\chi}_9\\
 n^9_{1,3} =2 \tilde{\chi} + 4 \tilde{\chi}_3 - 2 \tilde{\chi}_4 -2\tilde{\chi}_5 - \tilde{\chi}_8 +\tilde{\chi}_9 + n^0_{1,5}  - \eta\\
n^9_{2,7} = \eta\\
n^9_{1,4} = \tilde{\chi} + 3\tilde{\chi}_3 +  \tilde{\chi}_4 - 2\tilde{\chi}_5+\tilde{\chi}_7-\tilde{\chi}_8 -\tilde{\chi}_9+\tilde{\chi}_{10} - 2\eta+\zeta\\
n^9_{2,9}= -2\tilde{\chi}_3+\tilde{\chi}_4+\tilde{\chi}_5-\tilde{\chi}_7+\tilde{\chi}_8+\tilde{\chi}_9-\tilde{\chi}_{10} -n^0_{1,5}+\eta - \zeta\\
n^9_{1,5}=2\tilde{\chi}_3-\tilde{\chi}_4-
\tilde{\chi}_5+\tilde{\chi}_7- \tilde{\chi}_8 -\tilde{\chi}_9+
\tilde{\chi}_{10}+2n^0_{1,5}-\eta+\zeta
\end{array} \right.$$

One has $$ \epsilon_{10} =
-\tilde{\chi}_3+2\tilde{\chi}_6+\tilde{\chi}_{10}-\tilde{\chi}_{11}-n^0_{1,5}
-\eta$$ and then $B^{(10)}$ consists of following
coefficients:
$$B^{(10)} \left\{
\begin{array}{l}
n^{10}_{1,2} =2 \tilde{\chi}  - 2 \tilde{\chi}_3 +\tilde{\chi}_4  - 2 \tilde{\chi}_5 - \tilde{\chi}_7+ \tilde{\chi}_8+n^0_{1,5} + \eta-\zeta  \\
n^{10}_{4,9}=\zeta \\
 n^{10}_{3,7} =\tilde{\chi}  -\tilde{\chi}_3 +  \tilde{\chi}_6 + \tilde{\chi}_7- \tilde{\chi}_8  - \eta-\zeta\\
 n^{10}_{2,5} = \tilde{\chi} +\tilde{\chi}_3+\tilde{\chi}_4+  \tilde{\chi}_5 - 3 \tilde{\chi}_6 - \tilde{\chi}_7+ \tilde{\chi}_9- n^0_{1,5} + \eta \\
 n^{10}_{3,8}=-\tilde{\chi}_3-\tilde{\chi}_4+\tilde{\chi}_5+\tilde{\chi}_6+\tilde{\chi}_8-\tilde{\chi}_9\\
 n^{10}_{1,3} =2 \tilde{\chi} + 5 \tilde{\chi}_3 - 2 \tilde{\chi}_4 -2\tilde{\chi}_5 -2\tilde{\chi}_6- \tilde{\chi}_8 +\tilde{\chi}_9 -\tilde{\chi}_{10}+\tilde{\chi}_{11}+ 2n^0_{1,5} \\
 n^{10}_{3,10}=-\tilde{\chi}_3+2\tilde{\chi}_6+\tilde{\chi}_{10}-\tilde{\chi}_{11} -n^0_{1,5}-\eta \\
n^{10}_{2,7} = \tilde{\chi}_3-2\tilde{\chi}_6-\tilde{\chi}_{10}+\tilde{\chi}_{11} +n^0_{1,5}+2\eta\\
n^{10}_{1,4} = \tilde{\chi} + 3\tilde{\chi}_3 +\tilde{\chi}_4 -2\tilde{\chi}_5 +\tilde{\chi}_7-\tilde{\chi}_8-\tilde{\chi}_9+\tilde{\chi}_{10} -2\eta + \zeta\\
n^{10}_{2,9}= -2\tilde{\chi}_3+\tilde{\chi}_4+\tilde{\chi}_5-\tilde{\chi}_7+\tilde{\chi}_8+\tilde{\chi}_9-\tilde{\chi}_{10} -n^0_{1,5}+\eta - \zeta\\
n^{10}_{1,5}=2\tilde{\chi}_3-\tilde{\chi}_4-
\tilde{\chi}_5+\tilde{\chi}_7- \tilde{\chi}_8 -\tilde{\chi}_9+
\tilde{\chi}_{10}+2n^0_{1,5}-\eta+\zeta
\end{array} \right.$$

By computing $\Delta^{11}(B^{(10)})$, we get
$$ \epsilon_{11}= \tilde{\chi}-\tilde{\chi}_3+\tilde{\chi}_4-\tilde{\chi}_7+\tilde{\chi}_9+\tilde{\chi}_{11}-\tilde{\chi}_{12}-n^0_{1,5}-\zeta.$$
Let $\alpha$ be the number of prime packings of type
$\{(1,2),(4,9)\} \succ \{(5,11)\}$ and $\beta$ be the number of
prime packings of type $\{(1,3),(3,8)\} \succ \{(4,11)\}$. Then we
get $B^{(11)}$ with{\footnotesize
$$B^{(11)}\left\{
\begin{array}{l}
 n^{11}_{1,2} =  2 \tilde{\chi}- 2 \tilde{\chi}_3+\tilde{\chi}_4- 2\tilde{\chi}_5 - \tilde{\chi}_7+ \tilde{\chi}_8+ n^0_{1,5} + \eta   - \zeta  - \alpha \\
  n^{11}_{5,11} = \alpha \\
  n^{11}_{4,9} = \zeta - \alpha\\
  n^{11}_{3,7} = \tilde{\chi}- \tilde{\chi}_3 +\tilde{\chi}_6+  \tilde{\chi}_7 -  \tilde{\chi}_8-  \eta - \zeta     \\
    n^{11}_{2,5} =  \tilde{\chi} + \tilde{\chi}_3+ \tilde{\chi}_4+ \tilde{\chi}_5 - 3 \tilde{\chi}_6-\tilde{\chi}_7  + \tilde{\chi}_9  -n^0_{1,5}+ \eta \\
  n^{11}_{3,8} = - \tilde{\chi}_3-\tilde{\chi}_4  + \tilde{\chi}_5 + \tilde{\chi}_6 + \tilde{\chi}_8- \tilde{\chi}_9  - \beta\\
  n^{11}_{4,11} = \beta\\
  n^{11}_{1,3} = 2 \tilde{\chi} +5 \tilde{\chi}_3 - 2 \tilde{\chi}_4-2 \tilde{\chi}_5- 2 \tilde{\chi}_6 - \tilde{\chi}_8 + \tilde{\chi}_9- \tilde{\chi}_{10}+ \tilde{\chi}_{11}+ 2 n^0_{1,5}-\beta\\
n^{11}_{3,10}= -\tilde{\chi}_3 + 2 \tilde{\chi}_6 + \tilde{\chi}_{10}- \tilde{\chi}_{11} - n^0_{1,5} - \eta \\
  n^{11}_{2,7} = - \tilde{\chi}+ 2 \tilde{\chi}_3- \tilde{\chi}_4- 2 \tilde{\chi}_6 + \tilde{\chi}_7-\tilde{\chi}_9- \tilde{\chi}_{10}+ \tilde{\chi}_{12}+ 2 n^0_{1,5}+ 2\eta+\zeta+ \alpha + \beta\\
  n^{11}_{3,11} = \tilde{\chi} - \tilde{\chi}_3 + \tilde{\chi}_4 - \tilde{\chi}_7 + \tilde{\chi}_9 + \tilde{\chi}_{11}- \tilde{\chi}_{12}- n^0_{1,5} - \zeta   - \alpha - \beta\\
  n^{11}_{1,4} = 4 \tilde{\chi}_3- 2 \tilde{\chi}_5+2\tilde{\chi}_7-\tilde{\chi}_8-2\tilde{\chi}_9+ \tilde{\chi}_{10}-\tilde{\chi}_{11}+ \tilde{\chi}_{12}+n^0_{1,5}- 2\eta + 2\zeta + \alpha + \beta\\
 n^{11}_{2,9}=  - 2 \tilde{\chi}_3 + \tilde{\chi}_4+ \tilde{\chi}_5- \tilde{\chi}_7 + \tilde{\chi}_8+ \tilde{\chi}_9- \tilde{\chi}_{10}- n^0_{1,5}+\eta -\zeta\\
n^{11}_{1,5}=2\tilde{\chi}_3-
\tilde{\chi}_4-\tilde{\chi}_5+\tilde{\chi}_7- \tilde{\chi}_8
-\tilde{\chi}_9+ \tilde{\chi}_{10}+2n^0_{1,5}-\eta+\zeta
\end{array}
\right. $$}

Finally since $$
  \epsilon_{12}= -\tilde{\chi}-3\tilde{\chi}_3+2\tilde{\chi}_5+\tilde{\chi}_6
  -\tilde{\chi}_7+\tilde{\chi}_8+\tilde{\chi}_{12}-\tilde{\chi}_{13}-
  n_{1,5}^0+\eta.$$ we get $B^{(12)}$ with{\footnotesize
$$B^{(12)}\left\{
\begin{array}{l}
n^{12}_{1,2} =  2 \tilde{\chi}- 2 \tilde{\chi}_3+\tilde{\chi}_4- 2\tilde{\chi}_5 - \tilde{\chi}_7+ \tilde{\chi}_8+ n^0_{1,5} + \eta   - \zeta  - \alpha \\
  n^{12}_{5,11} = \alpha \\
  n^{12}_{4,9} = \zeta - \alpha\\
  n^{12}_{3,7} = 2 \tilde{\chi}+2 \tilde{\chi}_3 - 2 \tilde{\chi}_5+ 2 \tilde{\chi}_7 - 2 \tilde{\chi}_8-\tilde{\chi}_{12}+ \tilde{\chi}_{13}- 2 \eta - \zeta  + n^0_{1,5}   \\
  n^{12}_{5,12} =  -\tilde{\chi} - 3 \tilde{\chi}_3+ 2 \tilde{\chi}_5  + \tilde{\chi}_6 - \tilde{\chi}_7 + \tilde{\chi}_8 + \tilde{\chi}_{12}-\tilde{\chi}_{13}+ \eta -n^0_{1,5} \\
  n^{12}_{2,5} = 2 \tilde{\chi} +4 \tilde{\chi}_3+ \tilde{\chi}_4- \tilde{\chi}_5 - 4 \tilde{\chi}_6 - \tilde{\chi}_8 + \tilde{\chi}_9 - \tilde{\chi}_{12} + \tilde{\chi}_{13} \\
  n^{12}_{3,8} = - \tilde{\chi}_3-\tilde{\chi}_4  + \tilde{\chi}_5 + \tilde{\chi}_6 + \tilde{\chi}_8- \tilde{\chi}_9  - \beta\\
  n^{12}_{4,11} = \beta\\
  n^{12}_{1,3} = 2 \tilde{\chi} +5 \tilde{\chi}_3 - 2 \tilde{\chi}_4-2 \tilde{\chi}_5- 2 \tilde{\chi}_6 - \tilde{\chi}_8 + \tilde{\chi}_9- \tilde{\chi}_{10}+ \tilde{\chi}_{11}+ 2 n^0_{1,5}-\beta\\
n^{12}_{3,10}= -\tilde{\chi}_3 + 2 \tilde{\chi}_6 + \tilde{\chi}_{10}- \tilde{\chi}_{11} - n^0_{1,5}- \eta  \\
  n^{12}_{2,7} = - \tilde{\chi}+ 2 \tilde{\chi}_3- \tilde{\chi}_4- 2 \tilde{\chi}_6 + \tilde{\chi}_7-\tilde{\chi}_9- \tilde{\chi}_{10}+ \tilde{\chi}_{12}+ 2 n^0_{1,5}+ 2\eta+\zeta+ \alpha + \beta\\
  n^{12}_{3,11} = \tilde{\chi} - \tilde{\chi}_3 + \tilde{\chi}_4 - \tilde{\chi}_7 + \tilde{\chi}_9 + \tilde{\chi}_{11}- \tilde{\chi}_{12}- n^0_{1,5} - \zeta   - \alpha - \beta\\
  n^{12}_{1,4} = 4 \tilde{\chi}_3- 2 \tilde{\chi}_5+2\tilde{\chi}_7-\tilde{\chi}_8-2\tilde{\chi}_9+ \tilde{\chi}_{10}-\tilde{\chi}_{11}+ \tilde{\chi}_{12}+n^0_{1,5}- 2\eta + 2\zeta + \alpha + \beta\\
 n^{12}_{2,9}=  - 2 \tilde{\chi}_3 + \tilde{\chi}_4+ \tilde{\chi}_5- \tilde{\chi}_7 + \tilde{\chi}_8+ \tilde{\chi}_9- \tilde{\chi}_{10}- n^0_{1,5}+\eta -\zeta\\
n^{12}_{1,5}=2\tilde{\chi}_3- \tilde{\chi}_4-\tilde{\chi}_5+\tilde{\chi}_7-
\tilde{\chi}_8 -\tilde{\chi}_9+ \tilde{\chi}_{10}+2n^0_{1,5}-\eta+\zeta
\end{array}
\right. \eqno(3.15)$$}

To recall the meaning of several symbols, $\eta$ is the number of
prime packings of type $\{(1,3),(1,4)\} \succ \{(2,7)\}$, $\zeta$ is
the number of prime packings of type $\{(1,2),(3,7)\} \succ
\{(4,9)\}$, $\alpha$ is the number of prime packings of type
$\{(1,2),(4,9)\} \succ \{(5,11)\}$ and $\beta$ is the number of
prime packings of type $\{(1,3),(3,8)\} \succ \{(4,11)\}$.

\section{\bf Main results on general type 3-folds}

In this section, we would like to utilize those equalities and
inequalities of formal baskets to study 3-folds of general type. Let
$V$ be a nonsingular projective 3-fold of general type. The
3-dimensional Minimal Model Program (cf. \cite{KMM, K-M, Reid83})
says that $V$ has a minimal model $X$ with $\bQ$-factorial terminal
singularities. Therefore to study the birational geometry of $V$ is
equivalent to study that of $X$.

Let us begin with recalling some known relevant results. The
following theorem was proved by the first author and Hacon.

\begin{thm}\label{q>0}  (\cite{JC-H}) Assume $q(X):=h^1(\OO_X)>0$. Then $P_m >0$ for all $m \ge
2$ and $\varphi_m$ is birational for all $m\geq 7$.
\end{thm}

Thus we do not need to worry about irregular 3-folds in the
following discussion. The following result is due to  Koll\'ar.

\begin{thm}\label{Kollar} (\cite[Corollary 4.8]{Kol})
Assume $P_{m_0}:=P_{m_0}(X)\geq 2$ for some integer $m_0>0$. Then
$\varphi_{11m_0+5}$ is birational onto its image.
\end{thm}

Koll\'ar's result was improved by the second author.

\begin{thm}\label{5k+6} (\cite[Theorem 0.1]{JPAA}) Assume $P_{m_0}:=P_{m_0}(X)\geq 2$
for some integer $m_0>0$. Then $\varphi_{m}$ is birational onto its image for all $m\geq
5m_0+6$.
\end{thm}

\begin{setup}{\bf Other known results.}\label{other}
\begin{itemize}
\item[(i)] When $X$ is Gorenstein, it is proved in \cite{Crelle} that $\varphi_m$ is birational
for all $m\geq 5$.

\item[(ii)]
 When $\chi(\OO_X)<0$, Reid's
formula (4.1) says $P_2\geq 4$ and $P_m
> 0$ for all $m \ge 2$. It is proved in \cite[Corollary 1.3]{Chen-Zuo} that $\varphi_m$ is birational for all $m\geq
8$.

\item[(iii)]
When $\chi(\OO_X)=0$, since one can verify $l_Q(3)\geq l_Q(2)$ for
any basket $Q$, Reid's formula (4.1) says: $P_3(X)>P_2(X)>0$.
Moreover, $P_{m+1} \ge P_m$ for all $m \ge 2$.  So $P_3(X)\geq 2$.
It is proved in \cite[Theorem 1.4]{Chen-Zuo} that $\varphi_m$ is
birational for all $m\geq 14$.
\end{itemize}
\end{setup}



\begin{setup} {}From now on, we only  study  minimal 3-fold $X$ of general type
with
$\chi(\OO_X)>0$. Recall that $X$ is always attached the formal basket ${\bf B}(X)$.
Moreover, since $X$ is minimal and of general type, the vanishing
theorem (\cite{Ka,V}) on $X$ gives $\chi_m(X)=P_m(X)$.
Therefore we have various equalities and inequalities among plurigenera by the results in the previous sections.
Furthermore,
the canonical volume $\Vol(V)=\Vol(X)$ is nothing but $K^3_X$.
\end{setup}

The following result is due to Iano-Fletcher.

\begin{thm}\label{chi=1}(\cite{Flt}) Assume $\chi(\OO_X)=1$. Then $P_{12} \ge
1$ and $P_{24} \ge 2$.
\end{thm}

Combining all known results, we only need to consider the 3-fold $X$
satisfying $\chi(\OO_X)\ge 2$ and  $P_m \leq 1$ for all $2 \le m \leq 12$.

\begin{thm}\label{finiteness}
There are only finitely many formal baskets of minimal threefolds of general type satisfying  $\chi \ge 2$ and  $P_m \leq 1$ for all $2 \le m \leq 12$.
\end{thm}

\begin{proof}
By looking at inequality $(3.14)$, we have $$8 \ge \chi(\OO_X)+R\geq
\chi(\OO_X)$$ since $1 \ge \chi_m(X)=P_m(X) \ge 0$ for all $2 \le m
\le 12$.  Moreover, $8 \ge R$ implies that $n^0_{1,r}=0$ for all $r
\ge 9$. By equality $(3.5)$, one has $\sigma=\sum_{r=2}^8
n^0_{1,r}=10 \chi+5 P_2-P_3 \le 85$. It's clear that there are
finitely many initial basket $B^0=\{ n^0_{1,r}\}$ satisfying $\sigma
\le 85$ and $n^0_{1,r}=0$ for all $r \ge 9$. Each initial basket
allows finite ways of packings. Hence it follows that there are only
finitely many formal baskets satisfying the given conditions.
\end{proof}

By Theorem \ref{finiteness}, one can  obtain various
effective results by working out the classification of formal
baskets with small plurigenera. Indeed, by some more careful usage
of those inequalities in the previous section, we are able to obtain
our main results without too much extra works.

\begin{lem} \label{p20}
If $P_m \leq 1$ for all $m \leq 12$, then $P_2=0$.
\end{lem}

\begin{proof}
Recalling equation (3.10), we have:  $$
\epsilon_6=-3P_2-P_3+P_4+P_5+P_6-P_7-\epsilon=0
$$
which is equivalent to
$$ P_4+P_5+P_6 = 3 P_2+P_3+P_7+\epsilon. $$

If $P_2=1$, then $P_4=P_5=P_6=1$. It follows that $P_3=P_7=\epsilon=0$. But this is impossible since
$P_2=P_5=1$ implies $P_7 \ge 1$.
\end{proof}




\begin{lem}\label{chi}
Assume that $\chi(\OO_X) \geq 2$ and $P_{m} \le  1$ for  $m \leq 12$. Then  $\chi(\OO_X) \leq 6$.
\end{lem}

\begin{proof}
If $P_m \le 1$ for all $m \le 12$, we have seen $P_2=0$. Then by
inequality $(3.14)$, we get $8 \ge \chi=\chi(\OO_X)$. If $\chi =7 $
or $8$, then $P_5=P_6=1$. It follows that $P_{10}=P_{11}=P_{12}=1$.
Hence $8 \ge \chi +1$ gives $\chi =7$ and $P_8=1$ as well. Then
$P_{13}=1$. This leads to $8 \ge \chi+2 =9$, a contradiction.
\end{proof}

\begin{thm} \label{p12} Let $X$ be a projective minimal
3-fold of general type. Then $P_{12}\geq 1$.
\end{thm}

\begin{proof} It suffices to prove this for the situation $\chi \ge 2$ by \ref{other}(ii), (iii) and
Theorem \ref{chi=1}. We assume $P_{12}=0$ and will deduce a
contradiction. It's then clear that $P_2=P_3=P_4=P_6=0$.
\medskip

\noindent {\bf Step 1.} If $P_5=0$, then the equality (3.10) for
$\epsilon_6$ gives $P_7=\epsilon=0$. This also means $\sigma_5=0$.
Hence Assumption \ref{r6} is satisfied. Now since $\epsilon_7 \geq
\eta$ and $\epsilon_{12} \ge 0$ (cf. $(3.11),(3.13)$), one gets
$$ \chi \ge P_8+\eta \ge \chi+P_{13}.$$
It follows that $\chi=P_8+\eta$, $\epsilon_7=\eta$ and
$n^7_{3,7}=0$. Since $n^9_{3,7}= -\zeta$, we have $\zeta=0$. Now
$n_{4,9}^{11}=\zeta-\alpha\geq 0$ gives $\alpha=0$.

Hence since $n_{1,5}^0=0$ and so $n^9_{2,9}=-n^9_{1,5}=0$, we have
$n^9_{2,9}=0$ and $\epsilon_9=n^9_{2,9}+\zeta=0$ which gives
$P_{10}=P_8+P_9+\eta$.

Now $n^{12}_{3,8}+n^{12}_{2,7} \ge 0$ gives $\eta \geq \chi+3P_9 =
\eta+P_8+3P_9$.
 Hence $P_8=P_9=0$, and also $P_{10}=\eta=\chi$.
However, $n^{12}_{3,8}+n^{12}_{1,4}=P_{10}-2\eta
-P_{11}=-\chi-P_{11}<0$, which is a contradiction.
\medskip

\noindent {\bf Step 2.} If $P_5 >0$, then we have $P_7=0$. First of
all, (3.10) gives $P_5=\epsilon:=n^0_{1,5}+2 \sum_{r \ge 6}
n^0_{1,r}$. Since $n_{1,4}^7\geq 0$, one has
$$\chi\geq P_5+\eta+\sigma_5.$$
Again $\epsilon_{12} \ge 0$ (cf. $(3.13)$) gives the inequality:
 $$ 2P_5+P_8+\eta
\ge \chi +P_{13} +(8
\sigma_5-7n^0_{1,5}-5n^0_{1,6}-5n^0_{1,7}-\cdots-n^0_{1,11}).$$
Combining these two inequalities, we get $$2 \epsilon +P_8+\eta =
2P_5+P_8+\eta\geq P_5 +P_{13}+\eta + R', $$ where $R'=
2n^0_{1,5}+4n^0_{1,6}+4n^0_{1,7}+\cdots+8 n^0_{1,11}+8 \sum_{r \ge
12} n^0_{1,r}\geq 2\epsilon$. It follows that $P_8\geq P_5+P_{13}$.
Since $P_5>0$, we get $P_8>0$ and thus $P_{13}\geq P_{8}$. This
means $P_5=0$, a contradiction.
\end{proof}

\begin{lem}\label{beau} Let $W$ be a projective variety with at worst canonical singularities.
Given positive integers $m$ and $n$. Let $l:=\text{l.c.m}(m,n)$ and
$d:=\text{g.c.d}(m,n)$. Suppose that $P_m=P_n=P_l=1$. Then $P_d=1$.
\end{lem}

\begin{proof} Let $\pi: \tilde{W} \to W$ be a resolution of singularities.
It's clear that $P_k(\tilde{W})=P_k(W)$ for all $k \ge 1$. We may
thus assume that $W$ is nonsingular. The same argument as in
\cite[Lemma VIII.1.c]{beau} concludes the statement.
\end{proof}


\begin{thm}\label{p24}
Let $X$ be a  projective minimal 3-fold of general type. Then either
$P_{10} \ge 2$ or $P_{24}\geq 2$.
\end{thm}

\begin{proof} By  \ref{other} and Theorem \ref{chi=1}, we may only study those 3-folds with $\chi=\chi(\OO_X) \ge 2$.
Suppose, on the contrary, that $P_{24}\le 1$ and $P_{10} \le 1$. By
Theorem \ref{p12}, one has $P_{12}=P_{24}=1$. We will deduce a
contradiction.


\medskip
\noindent
{\bf Claim 1.} If $P_{8}>0$, then $P_4=P_8=1$.

In fact, this follows from Lemma \ref{beau} by taking $m=12$ and
$n=8$.

\medskip
Set $d:=\min \{ m | P_m(X)>0, m\in \mathbb{Z}^+\}.$ Clearly, one has
$d\leq 12$.

\medskip
\noindent
{\bf Claim 2.} If  $d|24$, then $P_n=0$ for any positive integer $n
\le \frac{24}{d}$ with $\text{g.c.d}(n,d)<d$.

To see this, suppose that $P_n >0$ for some $n \leq \frac{24}{d}$
with $d \nmid n$. Since $P_d>0$ and $d|24$, we see that
$1=P_{24}\geq P_{nd}\geq P_n$.  Thus, for $l:=\text{l.c.m}(n,d)$,
$P_l=1$. Now Lemma \ref{beau} gives $P_{(n,d)}=1$, contradicting to
minimality of $d$.

\medskip
\noindent
{\bf Claim 3.} We may assume that $d \ge 3$, i.e. $P_2=0$.

If $d=1$, then $P_m=1$ for all $m \le 12$.
But equality (3.10) gives  $\epsilon_6=-2-\epsilon=0$, a
contradiction.

 If $d=2$, then $P_4=P_6=1$ and Claim 2 tells
that $P_3=P_5=P_7=0$. Again equality (3.10) gives
$\epsilon_6=-1-\epsilon=0$, a contradiction.

\medskip
In what follows, we are going to apply those formulae in Section 4.
Recall, from equality (3.9), that $\epsilon:=n_{1,5}^0+2 \sum_{r \ge
6} n^0_{1,r}$. We will frequently use the following:
\medskip

{\bf Observation.} If $\epsilon+P_7=1$, then one of the following
situations occurs:
\begin{itemize}
\item[(A).]  $P_7=1$ and $n^0_{1,r}=0$ for all $r \ge 5$.
\item[(B).] $P_7=0$,  $n^0_{1,5}=1$ and $n^0_{1,r}=0$ for all $r \ge 6$.
\end{itemize}

Thus Assumption \ref{r6} is satisfied under both situations.
\medskip

Now we are ready for the proof, which is the case-by-case analysis
though it is slightly long.
\medskip

\noindent{\bf Case 1.} If $d=3$. Then, sine $P_9\leq P_{12}$, one
has $P_3=P_6=P_9=1$. By Claim 2, one gets $P_4=P_5=P_7=P_8=0$. Now
equality (3.10) gives $\epsilon_6=-\epsilon=0$. It follows that
$n^0_{1,r}=0$ for all $r \geq 5$ and hence Assumption \ref{r6} is
satisfied. But then, one will get $\epsilon_8=-1$, a
contradiction.
\medskip

\noindent{\bf Case 2.} If $d=4$, then $P_4=P_8=1$. One has
$P_5=P_6=0$ by Claim 2. Now equality (3.10) gives
$P_7+\epsilon=1$. Thus Assumption \ref{r6} is satisfied and so
$P_9=0$ by the inequality $\epsilon_8=-P_9\geq  0$. We discuss the
two cases in Observation:
\medskip

{\bf (2-A)} If $P_7=1$ and $\epsilon=0$, then we have $P_{11} \geq
P_7 \ge 1$. Now $\epsilon_{10} \ge 0$ yields $$ P_{10}\geq
P_{11}+n^0_{1,5}+\eta \geq 1.$$ This means, by our assumption on
$P_{10}$, that $P_{10}=1$ and $n_{1,5}^0=\eta=0$. So inequality
(4.1) gives
$$ 3=P_8+P_{10}+P_{12} \geq \chi+1+P_{11}+P_{13}+R \geq \chi+2,$$
contradicting to our assumption $\chi \geq 2$.
\medskip

{\bf (2-B)} If $P_7=0$ and $\epsilon=1$, then $n^0_{1,5}=1$.
Again, $\epsilon_{10} \ge 0$ gives
 $$ P_{10}\geq  P_{11}+n^0_{1,5}+\eta\geq P_{11}+\eta+1.$$ Thus $P_{10}=1$ and $P_{11}=\eta=0$. So inequality (3.14)
  yields
$$ 3=P_8+P_{10}+P_{12} \ge \chi+P_{13}+R \geq \chi+2,$$
 contradicting to the assumption $\chi \ge 2$.
\medskip

\noindent{\bf Case 3.} If $d=7$, then $P_2=\cdots=P_6=0$. But then
equality (3.10) gives $\epsilon_6=-P_7-\epsilon<0$, a
contradiction.
\medskip

\noindent{\bf Case 4.} If $d=8$, then, by Claim 1, $P_4=1$, a
contradiction.
\medskip

\noindent{\bf Case 5.} If $d=9$, then (3.10) gives $\epsilon=0$.
Hence Assumption \ref{r6} is satisfied. Now $\epsilon_8=-P_9<0$
yields a contradiction.
\medskip

\noindent{\bf Case 6.} If $d=10$, then, similarly, (3.10) gives
$\epsilon=0$ and thus Assumption \ref{r6} is satisfied. Now
$\epsilon_9\geq 0$ and $\epsilon_{10} \ge 0$ imply:
$$\eta\geq P_{10}\geq P_{11}+\eta.$$
It follows that $\eta=P_{10}=1$ and $P_{11}=0$. So inequality (3.14)
gives $ 2 \ge \chi+P_{13}$, which implies $P_{13}=0$ and $\chi=2$.
Now the direct computation shows $\epsilon_{12}=0$ and thus
$$B^{(11)}=B^{(12)}=\{ 5 \times (1,2), (3,7), 3 \times (2,5), 3 \times (1,3),
(3,11)\}.$$ But now we see that $B^{(12)}$ admits no non-trivial
prime packing of level $>12$.
This already says
$B^{(12)}=B^{(13)}=\cdots =B$. Therefore, there is only one  the formal basket
${\bf B}=(B, \chi(\OO_X), P_2)=(B,2,0)$ in this case. By the direct computation,
we see $P_{10}({\bf B})=0$ and $P_{24}({\bf B})=8$, a
contradiction.
\medskip

\noindent{\bf Case 7.} If $d=11$, then (3.10) gives $\epsilon=0$
and hence Assumption \ref{r6} is satisfied. But then
$\epsilon_{10}=-P_{11}-\eta <0$, a contradiction.
\medskip

\noindent{\bf Case 8.} If $d=12$, then similarly (3.10) gives
$\epsilon=0$ and hence Assumption \ref{r6} is satisfied. But then
inequality (3.14) yields $1=P_{12} \geq \chi+P_{13}\geq 2$, a
contradiction to the assumption $\chi \geq 2$.
\medskip

\noindent {\bf Notation.} In what follows, we will abuse the notation of a basket $B$ with its associated formal basket ${\bf B}=(B,\chi,\chi_2)=(B,\chi,0)$.
\medskip

\noindent{\bf Case 9.} If $d=6$, then $P_8=0$ by Claim 1. Since
$0<P_6\leq P_{18}\leq P_{24}=1$, we have $P_9=0,1$. Suppose
$P_9=1$, then Lemma \ref{beau} gives $P_3=1$, a contradiction to
$d=6$. Hence we have seen $P_9=0$. Now $\epsilon_6=0$ implies $P_7
+\epsilon=1$. Thus we get two situations as follows:
\medskip

{\bf (11-A)} $(P_7, \epsilon)=(0,1)$. Then
 $\epsilon_9\geq 0$ and $\epsilon_{10} \ge 0$ give
$$\eta+1\geq P_{10}+2\geq P_{11}+\eta+1.$$ In
particular, one has $P_{11}=0$ and $\eta=P_{10}+1$. Recall that
$P_{10}\leq 1$ by our assumption.
\medskip

{\bf (11-B)} $(P_7, \epsilon)=(1,0)$. Then
 $\epsilon_9\geq 0$ and $\epsilon_{10} \ge 0$ give
$$\eta+1\geq P_{10}+2\geq P_{11}+\eta.$$
In particular, one has $1 \ge P_{11}$ and $P_{10}+2 \ge \eta \ge
P_{10}+1$.
\medskip

The following table is the summary on the possible value of
$(P_7,P_{10},P_{11})$. Note, however, that all items should be
non-negative by our definition. {\tiny
$$ \begin{array}{c|cccccc}
(P_7,P_{10},P_{11}) & (0,0,0) & (0, 1,0) & (1,0,0) & (1,0,1) & (1,1,0) & (1,1,1) \\
\hline
\epsilon_7 & \chi+1 & \chi+1 & \chi+2 & \chi+2 & \chi+2 & \chi+2 \\
\epsilon_8 & 1 & 1 & 1 & 1 & 1 & 1\\
\epsilon_9 & -1+\eta & -2+\eta & -1+\eta & -1+\eta & -2+\eta & -2+\eta\\
\epsilon_{10} & 1-\eta &2-\eta & 2-\eta & 1 -\eta & 3-\eta & 2-\eta \\
\epsilon_{10}+\epsilon_{12} & 2-\chi-P_{13} &3-\chi-P_{13}&3-\chi-P_{13}&2-\chi-P_{13}&4-\chi-P_{13}&3-\chi-P_{13}
\end{array}
$$}
We are going to discuss it case by case.
\medskip

\noindent{\bf Subcase 9-I.} $(P_7,P_{10},P_{11})=(0,0,0)$.

The table tells $2 \ge \chi$ and $\eta=1$, hence $\chi=2$. But then
$n^8_{2,5}=-1$, a contradiction.
\medskip

\noindent{\bf Subcase 9-II.} $(P_7,P_{10},P_{11})=(0,1,0)$.

The table tells $\eta=2$ and $3 \ge \chi$. If $\chi=2$, then
$n^7_{1,4}=-1$, a contradiction. Hence $\chi=3$. Then we see
$\epsilon_{12}=-P_{13}$, which means $P_{13}=0$ and thus
$\epsilon_{12}=0$. Also $n_{2,9}^{11}=-\zeta$ says $\zeta=0$. Then
$n_{4,9}^{11}=\zeta-\alpha\geq 0$ gives $\alpha=0$. Since
$n_{1,4}^{11}=\beta-1\geq 0$ and $n_{3,8}^{11}=1-\beta\geq 0$, we
have $\beta=1$. Now we have,
$$B^{(12)}=B^{(11)}=
\{9 \times (1,2), 2 \times (3,7), (2,5), (4,11), 4 \times (1,3), 2
\times (2,7), (1,5)\}.$$ The only 1-step prime packing of level $>12$  of $B^{(12)}$
can only happen between $(4,11)$ and $(1,3)$. We obtained
$$\hat{B}=\{9 \times (1,2), 2 \times (3,7), (2,5), (5,14), 3\times (1,3), 2
\times (2,7), (1,5)\}.$$  We see $K^3(\hat{{ B}})=0$, and thus $0>K^3({{ B}'})$ for any $\hat{B} \succ B'$ by Lemma \ref{inequal}.
Therefore,  we get
$B=B^{(12)}$. Thus $P_{24}(X)=P_{24}(B^{(12)})=6$, a contradiction.
\medskip

\noindent{\bf Subcase 9-III.} $(P_7,P_{10},P_{11})=(1,0,0)$.

We have $P_{13}\geq P_7 \ge 1$ since $P_6>0$. Thus the table tells
that $\eta=1,2$ and that $\chi \le 2$, hence $\chi=2$.

If $\eta=1$, then $n^8_{2,5}=-1$, a contradiction. If $\eta=2$, then
$\epsilon_9=1$. Since $n^9_{1,4}=-1+\zeta \ge 0$ while $\epsilon_9
\geq \zeta$, one sees that $\zeta=1$. It follows that
$\epsilon_{11}=-1<0$, a contradiction.
\medskip

\noindent{\bf Subcase 9-IV.} $(P_7,P_{10},P_{11})=(1,0,1)$.

Since $P_{13} \ge P_7\geq 1$, the table gives $\chi \le 1$, a
contradiction to $\chi\geq 2$.
\medskip

\noindent{\bf Subcase 9-V.} $(P_7,P_{10},P_{11})=(1,1,0)$.

Since $P_{13}>0$, the table tells that $\chi\leq 3$ and $2\leq
\eta\leq 3$.

If $\chi=2$ and $\eta=3$, then $n^7_{1,4}=2-\eta=-1<0$, a
contradiction.

If $\chi=3$ and $\eta=2$, then $\epsilon_{10}=1$ and
$\epsilon_{10}+\epsilon_{12}=0$. Thus $\epsilon_{12}=-1<0$, a
contradiction.

If $\chi=\eta=2$, we can determine other unknown quantities. First,
$n_{2,5}^{12}=-1+P_{13}\geq 0$ gives $P_{13}=1$. Thus
$\epsilon_{12}=0$ and $B^{(12)}=B^{(11)}$. Now
$n_{2,9}^{11}=-\zeta\geq 0$ gives $\zeta=0$. Then $n_{4,9}^{11}\geq
0$ tells $\alpha=0$. Finally $n_{3,11}^{11}=-\beta\geq 0$ says
$\beta=0$. Hence we get:
$$B^{(12)}=\{5 \times (1,2), 2 \times (3,7), (3,8),  (1,3), (3,10),
(2,7)\}.$$ It is clear that $B^{(12)}$ admits two 1-step prime
packings of level $>12$:
$$B'=\{5 \times (1,2), 2 \times (3,7), (3,8),  (1,3), (5,17)\}.$$
$$B''=\{5 \times (1,2), 2 \times (3,7), (3,8),  (4,13), (2,7)\}.$$
But $K^3(B'')<0$, $K^3(B')>0$ and $B'$ is a minimal positive formal basket, we see that either  $B^{(12)}\succcurlyeq B\succcurlyeq
B'$ or $B^{(12)}=B$. By a direction calculation, we see
$P_{24}(B^{(12)})=4$ and $P_{24}(B')=3$. Thus Lemma \ref{inequal}
implies $P_{24}=P_{24}(X)\geq 3$, a contradiction.

If $\chi=\eta=3$, then the table shows
$\epsilon_{10}=\epsilon_{12}=0$ and $P_{13}=1$. We detect $B^{(11)}$
as before. First, $n_{2,9}^{11}\geq 0$ and $n_{1,5}^{11}\geq 0$
imply $\zeta=1$. Then $\epsilon_{11}=1-\zeta=0$ implies
$\alpha=\beta=0$. So we get:
$$B^{(12)}=B^{(11)}=\{7 \times (1,2),(4,9),(3,7), 2 \times (2,5), (3,8),3
\times  (1,3), 3 \times (2,7)\}.$$ We see that $B^{(12)}$ admits
only two 1-step prime packings of level $>12$:
$$\hat{B}'=\{7 \times
(1,2),(7,16), 2 \times (2,5), (3,8), 3 \times  (1,3), 3 \times
(2,7)\},$$
$$\hat{B}''=\{7 \times (1,2), (4,9),( 3,7), (2,5), (5,13), 3 \times  (1,3), 3 \times (2,7)\}.$$
By computation, both $\hat{B}'$ and $\hat{B}''$ are minimal positive
(with regard to $B^{(12)}$). So we see that either
$B^{(12)}\succcurlyeq B\succcurlyeq \hat{B}'$ or
$B^{(12)}\succcurlyeq B\succcurlyeq \hat{B}''$. Since
$P_{24}(B^{(12)})=8$, $P_{24}(\hat{B}')=6$ and
$P_{24}(\hat{B}'')=4$, Lemma \ref{inequal} implies $P_{24}\geq 4$, a
contradiction.
\medskip

\noindent{\bf Subcase 9-VI.} $(P_7, P_{10},P_{11})=(1,1,1)$.

Since $P_{13}>0$, the table tells $\chi=2$, $\eta=2$ and
$\epsilon_{12}=0$. Now $n_{2,9}^{11}=-\zeta\geq 0$ gives $\zeta=0$.
Further, $n_{4,9}^{11}\geq 0$ gives $\alpha=0$. Finally,
$n_{3,8}^{11}=1-\beta\geq 0$ and $n_{1,4}^{11}=-1+\beta\geq 0$
implies $\beta=1$. So we have seen:
$$B^{(12)}=B^{(11)}=\{5 \times (1,2), 2 \times (3,7), (4,11),  (1,3), 2
\times (2,7)\}.$$ The only prime packing of $B^{(12)}$ of level $>12$ is the
following basket: $$B':=\{5 \times (1,2), 2 \times (3,7), (5,14), 2
\times (2,7)\}$$ with $K^3(B')=0$. This means that $B^{(12)}$ is
already minimal positive and thus $B=B^{(12)}$. So
$P_{24}=P_{24}(B^{(12)})=6>1$, a contradiction.
\medskip

\noindent{\bf Case 10.} If $d=5$, then Claim 1 implies $P_8=0$.
Also, we have $P_5\leq P_{10}\leq 1$, which means $P_5=1$.

{}First we study $P_6$. Assume $P_6>0$, then $P_6=1$ since $P_6\leq
P_{12}$. Since $0<P_6\leq P_{18}\leq P_{24}=1$, we have $P_9=0$,
$1$. Suppose $P_9=1$, then Lemma \ref{beau} gives $P_3=1$, a
contradiction to $d=5$. Hence we have seen $P_9=0$. Similarly, if
$P_8>0$, then $P_8=1$ since $P_8\leq P_{24}$. Lemma \ref{beau} gives
$P_2=1$, a contradiction to $d=5$. Thus $P_8=0$. Noting that
$P_{11}\geq P_6=1$, the inequality $\epsilon_9+\epsilon_{10} \ge 0$
gives:
$$P_5+1\geq P_7+9\sigma_5-(7n^0_{1,5}+6n^0_{1,6}+5n^0_{1,7}+3n^0_{1,8}+n^0_{1,9}). \eqno{(4.3)}$$
On the other hand, equality (3.10) implies:
$$P_5+1=P_7+\epsilon=P_7+\sigma_5+ \sum_{r \ge 6} n^0_{1,6}.\eqno{(4.4)}$$
Now (4.3) and (4.4) imply $n^0_{1,r}=0$ for all $r \ge 5$ and
$P_7=P_5+1 \ge 2$. It follows that $P_{12} \ge 2$, a contradiction.
Therefore we have actually shown that $P_6=0$.

Next we study $P_7$. Clearly $P_7\leq P_{12}=1$. Assume $P_7=0$.
Then equality (3.10) gives $\epsilon=1$. This corresponds to
Observation (B).  Now $\epsilon_9+\epsilon_{10} \ge 0$ implies that
$$1+P_9=P_5+P_9 \geq P_{11}+2.$$
Since $P_{15}>0$, we see $P_9\leq P_{24}=1$. Hence $P_9=1$,which
implies $P_{11}=0$. Now $\epsilon_{10}=-\eta$ gives $\eta=0$. Thus
we can see $\epsilon_9=0$. It follows that $\zeta=0$ since $\zeta
\le \epsilon_9$. Finally we can see that
$n^{11}_{2,7}+n^{11}_{4,9}+n^{11}_{3,8}=-\chi+1\leq -1$, which is a
contradiction. We have shown $P_7=1$.

To make a summary, we have: $P_5=P_7=P_{10}=P_{12}=1$ and
$P_2=P_3=P_4=P_6=P_8=0$. Note also that (3.10) gives $\epsilon=0$,
thus Assumption \ref{r6} is always satisfied. We need to study
$P_9$, $P_{11}$.

Clearly, $P_9\leq P_{24}=1$ since $P_{15}>0$. Again,
$\epsilon_9+\epsilon_{10} \geq 0$ gives $P_9 \ge P_{11}$. The next
table is a summary for 3 possibilities of $(P_9, P_{11})$:
$$ \begin{array}{c|ccc}
(P_9,P_{11}) & (0,0) & (1, 0) & (1,1)  \\
\hline
\epsilon_7 & \chi+1 & \chi+1 & \chi+1  \\
\epsilon_8 & 1 & 0 & 0 \\
\epsilon_9 & -1+\eta & \eta & +\eta \\
\epsilon_{10} & 1-\eta &1-\eta & -\eta  \\
\epsilon_{10}+\epsilon_{12} & 3-\chi-P_{13} &3-\chi-P_{13}&2-\chi-P_{13}
\end{array}
$$
\medskip

\noindent{\bf Subcase 10-I.} $(P_9,P_{11})=(0,0)$.

The table tells us that $\eta=1$ and $\chi=2$, $3$.

When $\chi=2$, $\epsilon_{11}=-\zeta\geq 0$ gives $\zeta=0$ and
thus $\epsilon_{11}=0$. This says $\alpha=\beta=0$. Since
$P_{13}\leq 1$ by the table, we first assume $P_{13}=0$. Then we
get
$$B^{(12)}=\{2\times (1,2), (3,7), (5,12), 2\times (2,5), (3,8), (1,3), (2,7) \}.$$
But we see that $K^3(B^{(12)})<0$, contradicting to
$K^3(B^{(12)})\geq K^3(B)=K_X^3>0$. Thus $P_{13}=1$,
$\epsilon_{12}=0$ and we get
$$B^{(12)}=\{2 \times (1,2), 2 \times (3,7) , 3 \times (2,5), (3,8), (1,3), (2,7) \}.$$
Since any further prime packing dominated by $B^{(12)}$ has
negative volume (due to the direct computation) and
$B^{(12)}\succcurlyeq B$, we see $B=B^{(12)}$. So
$P_{24}=P_{24}(B^{(12)})=4>1$, a contradiction.

When $\chi=3$, the table tells $P_{13}=0$ and $\epsilon_{12}=0$.
Since $n_{2,9}^{11}=-\zeta\geq 0$, we have $\zeta=0$. Thus by
$n_{4,9}^{11}=\zeta-\alpha\geq 0$, we see $\alpha=0$. Finally
$\epsilon_{11}=1$ gives $\beta\leq 1$. If $\beta=1$, then we get:
$$B^{(12)}=\{4 \times (1,2), 3 \times (3,7) , 4 \times (2,5), (4,11), 2 \times (1,3), (2,7), (1,4)\}.$$
But we see that $K^3(B^{(12)})<0$, contradicting to
$K^3(B^{(12)})\geq K^3(B)=K_X^3>0$. Thus we must have $\beta=0$.
Then we get:
$$B^{(12)}=\{4 \times (1,2), 3 \times (3,7) , 4 \times (2,5), (3,8), 3 \times (1,3), (3,11) \}.$$
Since any further prime packing dominated by $B^{(12)}$ has
negative volume (due to the direct computation) and
$B^{(12)}\succcurlyeq B$, we see $B=B^{(12)}$. So
$P_{24}=P_{24}(B^{(12)})=2>1$, a contradiction.
\medskip

\noindent{\bf Subcase 10-II.} $(P_9,P_{11})=(1,0)$.

The table tells that $\eta=0$, $1$ and $\chi=2$, $3$.

If $\eta=0$, then $n^{10}_{2,7}=-1$, a contradiction.


If $(\eta,\chi)=(1,2)$, then $n_{3,8}^{11}=-\beta\geq 0$ gives
$\beta=0$. Furthermore, $n_{4,9}^{11}+n_{3,11}^{11}=1-2\alpha\geq
0$ implies $\alpha=0$. Also, $n_{2,9}^{11}=1-\zeta\geq 0$ and
$n_{1,5}^{11}=\zeta-1\geq 0$ imply $\zeta=1$. Finally, the table
tells $\epsilon_{10}+\epsilon_{12}=1-P_{13}$ and so $P_{13}\leq 1$.
When $P_{13}=1$, we get:
$$B^{(12)}=\{ (1,2), (4,9), (3,7), 4 \times (2,5), 2 \times (1,3), (2,7) \}.$$
Clearly, $B^{(12)}$ admits only one prime packing of level $>12$:
$$\tilde{B}=\{ (1,2), (7,16), 4 \times (2,5), 2 \times (1,3), (2,7) \}.$$
Thus we see that either $B^{(12)}=B$ or $B^{(12)}\succcurlyeq
B\succcurlyeq \tilde{B}$. By computation, we know
$P_{24}(B^{(12)})=5$ and $P_{24}(\tilde{B})=3$. Thus
$P_{24}=P_{24}(B)\geq 3>1$, a contradiction. When $P_{13}=0$, we
get:
$$B^{(12)}=\{ (1,2), (4,9), (5,12), 3\times (2,5), 2 \times (1,3), (2,7) \}.$$
But we see that $K^3(B^{(12)})<0$, contradicting to
$K^3(B^{(12)})\geq K^3(B)=K_X^3>0$.


If $(\eta,\chi)=(1, 3)$, the table tells $P_{13}=0$ and
$\epsilon_{12}=0$. Also, $n_{2,9}^{11}=1-\zeta\geq 0$ and
$n_{1,5}^{11}=\zeta-1\geq 0$ imply $\zeta=1$. Furthermore,
$n_{3,8}^{11}=-\beta\geq 0$ gives $\beta=0$. Finally,
$n_{4,9}^{11}=1-\alpha\geq 0$ imply $\alpha\leq 1$. When
$\alpha=1$, we get:
$$B^{(12)}=\{ 2 \times (1,2), (5,11), 2 \times (3,7), 5 \times (2,5), 4 \times (1,3), (2,7), (1,4) \}.$$
But we see that $K^3(B^{(12)})<0$, contradicting to
$K^3(B^{(12)})\geq K^3(B)=K_X^3>0$. When $\alpha=0$, we get:
$$B^{(12)}=\{ 3 \times (1,2), (4,9), 2 \times (3,7), 5 \times (2,5), 4 \times (1,3), (3,11) \}.$$
There is  only one prime packing of level $>12$:
$$\{ 3 \times (1,2), (7,16), (3,7), 5 \times (2,5), 4 \times (1,3), (3,11) \},$$
which has $K^3<0$, we see that $B^{(12)}=B$. Thus
$P_{24}=P_{24}(B)=P_{24}(B^{(12)})=3>1$, a contradiction.
\medskip

\noindent{\bf Subcase 10-III.} $(P_9,P_{11})=(1,1)$.

The table tells that $\eta=0$ and $\chi =2$. Also $\epsilon_{10}=0$
implies $\zeta=0$. But then
$n^{11}_{2,7}+n^{11}_{4,9}+n^{11}_{3,8}=-2$, which is a
contradiction.
\medskip

This completes the proof.
\end{proof}

\begin{setup}{\bf Proof of Corollary \ref{birat}.}
\end{setup}
\begin{proof} (1) By virtue of \ref{other}, we may only study a
minimal 3-fold $X$ with $\chi(\OO_X)>0$. Then Theorem \ref{chi=1}
and Theorem \ref{p24} imply that there is a positive integer
$m_0\leq 24$ such that $P_{m_0}\geq 2$. Thus, by Theorem
\ref{5k+6}, $\varphi_m$ is birational for all $m\geq 126$.

(2) Set $\Phi:=\varphi_{126}$. By taking a proper birational
modification $\pi: \tilde{X}\rightarrow X$ (to resolve the
indeterminancy of $\Phi$), we may assume that $\Phi\circ
\pi:\tilde{X}\longrightarrow \bP^N$ is a birational morphism.
Denote by $M$ the movable part of $|126K_{\tilde{X}}|$. Then
$126\pi^*(K_X)\geq M=\Phi\circ \pi^*(H)$ for a very ample divisor
$H$ on $\bP^N$. Note that the image of $\tilde{X}\neq \bP^3$, we
see that $N>3$ and that:
$$(126\pi^*(K_X))^3\geq H^3\geq 2,$$
which at least gives $K_X^3\geq \frac{1}{63\cdot 126^2}$. We are
done.
\end{proof}

\begin{rem} We will develop some more methods and more detail classification to estimate the lower bound of $K_X^3$ in our next paper,
where a sharp bound is obtained.
To curb the length of this paper, we have to cut out other details
here.
\end{rem}



\begin{thebibliography}{99}
\bibitem{beau} A. Beauville, {\em Complex algebraic surfaces}. London
Mathematical Society Lecture Note Series, {\bf 68}. Cambridge
University Press, Cambridge, 1983.

\bibitem{Bom} E. Bombieri, {\em
Canonical models of surfaces of general type}. Inst. Hautes Etudes
Sci. Publ. Math. {\bf 42} (1973), 171--219.
\bibitem{Crelle} J. A. Chen, M. Chen, D.-Q. Zhang,
{\em The 5-canonical system on 3-folds of general type}. J. Reine
Angew. Math. {\bf 603} (2007), 165-181.

\bibitem{JC-H} J. A. Chen, C. D. Hacon, {\em Pluricanonical
systems on irregular 3-folds of general type}. Math. Z. {\bf
255}(2007), no. 2, 343-355

\bibitem{JLMS} J. A. Chen, M. Chen, {\em The
Canonical volume of 3-folds of general type with $\chi\leq 0$}. J.
London Math. Soc. {\bf (78)} 2008, 693-706.


\bibitem{IJM} M. Chen, {\em Canonical stability of 3-folds of general
type with $p_g\geq 3$}, Internat. J. Math. {\bf 14} (2003), 515-528

\bibitem{JPAA} M. Chen, {\em
On the ${\bf Q}$-divisor method and its application}. J. Pure Appl.
Algebra {\bf 191} (2004), 143-156.


\bibitem{MathAnn} M. Chen,  {\em A sharp lower bound for the canonical
 volume of 3-folds of general type}. Math. Annalen {\bf 337}(2007), no.4, 887-908

\bibitem{Chen-Zuo} M. Chen, K. Zuo, {\em Complex projective
threefolds with non-negative canonical Euler-Poincare
characteristic}, Comm. Anal. Geom. {\bf 16} (2008), 159-182.


\bibitem{Ein-Lazars} L. Ein, R. Lazarsfeld, {\em Global generation of
pluricanonical and adjoint linear series on smooth projective
threefolds}. J. Amer. Math. Soc. {\bf 6} (1993), no. 4, 875-903


\bibitem{H-M} C. D. Hacon and J. M$^{\text{\rm c}}$Kernan, {\em
Boundedness of pluricanonical maps of varieties of general type},
Invent. Math. 166 (2006), 1-25.

\bibitem{C-R} A. R. Iano-Fletcher,
{\em Working with wighted complete intersections}, Explicit
birational geometry of 3-folds. London Mathematical Society, Lecture
Note Series, 281. Cambridge University Press, Cambridge, 2000.

\bibitem{Fletcher} A. R. Iano-Fletcher,
{\em Inverting Reid's exact plurigenera formula}. Math. Ann. 284
(1989), no. 4, 617-629.

\bibitem{Flt}
A. R. Iano-Fletcher. {\em Contributions to Riemann-Roch on
Projective 3-folds with Only Canonical Singularities and
Applications}. Proceedings of Symposia in Pure Mathematics {\bf
46}(1987), 221-231.


\bibitem{Ka} Y. Kawamata, {\em A generalization of Kodaira-Ramanujam's
vanishing theorem}, Math. Ann. {\bf 261}(1982), 43-46.

\bibitem{Ka2} Y. Kawamata, {\em On the plurigenera of minimal algebraic $3$-folds with $K\equiv 0$},  Math. Ann.  {\bf 275}  (1986),  no. 4, 539--546.

\bibitem{KMM} Y. Kawamata, K. Matsuda, K. Matsuki, {\em Introduction to the
minimal model problem}, Adv. Stud. Pure Math. {\bf 10}(1987),
283-360.

\bibitem{Kol} J. Koll\'ar, {\em Higher direct images of dualizing
sheaves I}, Ann. Math. {\bf 123}(1986), 11-42;\ \ II, ibid. {\bf
124}(1986), 171-202.

\bibitem{K-M} J. Koll\'ar, S. Mori, Birational geometry of algebraic
varieties, 1998, Cambridge Univ. Press.

\bibitem{Luo} T. Luo, {\em Global $2$-forms on regular $3$-folds of general type}.
Duke Math. J. {\bf 71} (1993), no. 3, 859-869.

\bibitem{C-3f} M. Reid, {\em Canonical threefolds}, In A.
Beauville, editor, {\it G\'eom\'etric Alg\'ebrique Angers}, pages
273-310. Sijthoff \& Noordhoff, 1980.

\bibitem{Reid83} M. Reid, {\em Minimal models of
canonical 3-folds}, Adv. Stud. Pure Math. {\bf 1}(1983), 131-180.

\bibitem{YPG} M. Reid,  {\em Young person's guide to canonical
singularities}, Proc. Symposia in pure Math. {\bf 46}(1987),
345-414.

\bibitem{Ta} S. Takayama, {\em Pluricanonical systems on algebraic
varieties of general type}, Invent. Math. 165 (2006), 551 -- 587.

\bibitem{Tsuji} H. Tsuji, {\em Pluricanonical systems of projective
varieties of general type. I}. Osaka J. Math. {\bf 43} (2006), no.
4, 967-995
\bibitem{V} E. Viehweg, {\em Vanishing theorems}, J. reine angew. Math. {\bf
335}(1982), 1-8.

{\bf 113}(1985), 23-51.
\end{thebibliography}
\end{document}